\documentclass[12 pt]{article}
\usepackage{array,color}
\usepackage{amsbsy}
\usepackage{amsmath}
\usepackage{amssymb}
\usepackage{amsfonts}
\usepackage{latexsym}
\usepackage{epsf}
\usepackage{epsfig}
\usepackage{graphicx}
\usepackage{fancybox,fancyhdr,times,subfigure}
\usepackage{epsfig, color, graphicx, rotate, colordvi}
\usepackage{a4wide}
\usepackage[numbers]{natbib}
\usepackage{mathrsfs}
\usepackage{bm}

\usepackage{array,color}
\usepackage{amsbsy}
\usepackage{amsmath}
\usepackage{amssymb}
\usepackage{amsfonts}
\usepackage{latexsym}
\usepackage{epsf}
\usepackage{epsfig}
\usepackage{graphicx}
\usepackage{epstopdf}
\usepackage{fancybox,fancyhdr,times,subfigure}
\usepackage{epsfig, color, graphicx, rotate, colordvi}
\usepackage{a4wide}
\usepackage[numbers]{natbib}
\usepackage{mathrsfs}
\usepackage{bm}
\usepackage{float}
\usepackage{paralist}

\newcommand{\bbD}{{\bf D}}
\newcommand{\bbB}{{\bf B}}
\newcommand{\bbG}{{\bf G}}

\newcommand{\bbf}{{\bf f}}
\newcommand{\bbx}{{\bf x}}

\newcommand{\bbX}{{\bf X}}
\newcommand{\bby}{{\bf y}}

\newcommand{\bbI}{{\bf I}}
\newcommand{\bbJ}{{\bf J}}
\newcommand{\bbK}{{\bf K}}
\newcommand{\bbO}{{\bf O}}
\newcommand{\bbP}{{\bf P}}

\newcommand{\bbS}{{\bf S}}

\newcommand{\bbW}{{\bf W}}
\newcommand{\bbR}{{\bf R}}
\newcommand{\bbPhi}{{\bf \Phi}}

\newcommand{\bgL}{{\bf \Lambda}}
\newcommand{\bbF}{{\bf F}}
\newcommand{\bbH}{{\bf H}}

\newcommand{\bbe}{{\bf e}}
\newcommand{\bbb}{{\bf b}}
\newcommand{\bbQ}{{\bf Q}}
\newcommand{\bgma}{{\bf \gamma}}
\newcommand{\bbM}{{\bf M}}

\newcommand{\bqn}{\begin{eqnarray*}}
\newcommand{\eqn}{\end{eqnarray*}}
\newcommand{\rE}{{\textrm{E}}}
\newcommand{\rtr}{{\textrm{tr}}}
\newcommand{\bbA}{{\bf A}}
\newcommand{\bqa}{\begin{eqnarray}}
\newcommand{\eqa}{\end{eqnarray}}

\newcommand{\non}{\nonumber\\}
\newcommand{\ep}{\varepsilon}
\newcommand{\bbep}{\bf \varepsilon}
\newcommand{\rP}{\textrm{P}}
\newcommand{\ga}{\alpha}
\newcommand{\la}{\lambda}
\newtheorem{remark}{{\it Remark}}[section]
\newtheorem{lemma}{Lemma}[section]
\newcommand{\proof}{\textsc{Proof}}
\begin{document}
\title
{Order Determination of Large Dimensional Dynamic Factor Model}
\author{Z. D. Bai$^a$, Chen Wang$^b$, Ya Xue$^c$ and  Matthew Harding$^d$\\
Northeast Normal University$^a$, University of Cambridge$^b$\\
and University of California, Irvine$^{c,d}$\\[2ex]}
\maketitle
\begin{abstract}
Consider the following dynamic factor model:
$\bbR_{t}=\sum_{i=0}^{q}\bgL_{i}\bbf_{t-i}+\bbe_{t},t=1,...,T$,
where $\bgL_{i}$ is an $n\times k$ loading matrix of full rank, $\{\bbf_t\}$  are i.i.d. $k\times1$-factors, and $\bbe_t$ are independent $n\times1$ white noises. Now, assuming that $n/T\to c>0$, we want to estimate the orders $k$  and $q$ respectively. Define a random matrix
$$\bbPhi_{n}(\tau)=\frac{1}{2T}\sum_{j=1}^{T}(\bbR_{j}\bbR_{j+\tau}^{*}+\bbR_{j+\tau}\bbR_{j}^{*}),$$
where $\tau\ge 0$ is an integer. When there are no factors, the matrix $\Phi_{n}(\tau)$ reduces to
$$\bbM_n(\tau)=\frac{1}{2T}\sum_{j=1}^T(\bbe_{j}\bbe_{j+\tau}^{*}+\bbe_{j+\tau}\bbe_{j}^{*}).$$
When $\tau=0$, $\bbM_n(\tau)$ reduces to the usual sample covariance matrix whose ESD tends to the well known MP law and $\bbPhi_n(0)$ reduces to the standard spike model. Hence the number $k(q+1)$ can be estimated by the number of spiked eigenvalues of $\bbPhi_n(0)$. To obtain separate estimates of $k$  and $q$ , we have employed the spectral analysis of $\bbM_n(\tau)$ and established the spiked model analysis for $\bbPhi_n(\tau)$.
\end{abstract}
\section{Introduction}\label{sec1}
\setcounter{equation}{0}
\def\theequation{\thesection.\arabic{equation}}
\setcounter{subsection}{0}
For a $p\times p$ random Hermitian matrix $\bbA$ with eigenvalues $\la_j,j=1,2,\cdots, p$, the empirical spectral distribution (ESD) of $\bbA$ is defined as $$F^{\bbA}(x)=\frac{1}{p}\sum_{j=1}^pI(\la_j\leq x).$$ The limiting distribution $F$ of $\{F^{\bbA_n}\}$ for a given sequence of random matrices $\{\bbA_n\}$ is called the limiting spectral distribution (LSD).
Let $\{\ep_{it}\}$ be independent identically distributed (i.i.d) random variables with common mean 0, variance 1. Consider a high dimensional dynamic $k$-factor model with lag $q$, that is, $\bbR_{t}=\sum_{i=0}^{q}\bgL_{i}
\bbf_{t-i}+\bbe_{t},t=1,...,T$, where $\bgL_{i}$ is an $n\times k$ loading matrix of full rank, $\{\bbf_t\}$  are i.i.d. $k\times1$-factors with common mean 0, variance 1, whereas $\bbe_t$ corresponds to the noise component with $\bbe_t=(\ep_{1t},\cdots,\ep_{nt})'$. In addition, both components of $\bbe_t$ and $\bbf_t$ are assumed to have finite 4th moment.

This model can also be thought as an information-plus-noise type model (Dozier \& Silverstein, 2007a, b; Bai \& Silverstein, 2012). Here both $n$ and $T$ tend to $\infty$, with $n/T\to c$ for some $c>0$. Compared with $n$ and $T$, the number of factors $k$ and that of lags $q$ are fixed but unknown. An interesting and important problem to economists is how to estimate $k$ and $q$. To this end,  define $\bbPhi_{n}(\tau)=\frac{1}{2T}\sum_{j=1}^{T}(\bbR_{j}\bbR_{j+\tau}^{*}+\bbR_{j+\tau}\bbR_{j}^{*})$, $\bgma_t=\frac{1}{\sqrt{2T}}\bbe_t$ and $\bbM_n(\tau)=\sum_{k=1}^T(\bgma_k\bgma_{k+\tau}^*+\bgma_{k+\tau}\bgma_k^*), \tau=0,1,\cdots.$
Here $^*$ stands for the transpose and complex conjugate of a complex number and $\tau$ is referred to be the number of lags. Denote
\bqn
\boldsymbol{\bgL}&=&(\boldsymbol{\bgL}_0,\boldsymbol{\bgL}_1,\cdots,\boldsymbol{\bgL}_q)_{n\times k(q+1)},\\
\boldsymbol{\mathrm{F}^{\tau}}&=&
\left(
  \begin{array}{cccc}
    \bbf_{T+\tau} & \bbf_{T+\tau-1}& \cdots & \bbf_{\tau+1} \\
     \bbf_{T+\tau-1} & \bbf_{T+\tau-2}& \cdots & \bbf_{\tau} \\
    \vdots & \vdots & \vdots & \vdots \\
   \bbf_{T+\tau-q} & \bbf_{T+\tau-1-q}& \cdots & \bbf_{\tau+1-q} \\
  \end{array}
\right)
_{k(q+1)
\times T},\\
\boldsymbol{\bbe}^{\tau}&=&(\bbe_{T+\tau},\bbe_{T+\tau-1},\cdots,\bbe_{\tau+1})_{n\times T}.
\eqn
Then we have that
$\bbPhi_{n}(\tau)=\frac1{2T}[(\bgL\mathrm{F}^{\tau}+\boldsymbol{\bbe}^{\tau})(\bgL\mathrm{F}^{0}+\boldsymbol{\bbe}^{0})^*
+(\bgL\mathrm{F}^{0}+\boldsymbol{\bbe}^{0})(\bgL\mathrm{F}^{\tau}+\boldsymbol{\bbe}^{\tau})^*]$ and
$\bbM_n(\tau)=\frac1{2T}(\boldsymbol{\bbe}^{\tau}\boldsymbol{\bbe}^{0*}+\boldsymbol{\bbe}^{0}\boldsymbol{\bbe}^{\tau*})$.

Note that essentially, $\bbM_n(\tau)$ and $\Phi_{n}(\tau)$ are symmetrized auto-cross covariance matrices at lag $\tau$ and generalize the standard sample covariance matrices $\bbM_{n}(0)$ and $\Phi_{n}(0)$, respectively. The matrix $\bbM_{n}(0)$ has been intensively studied in the literature and it is well known that the LSD has an MP law (Mar\v{c}enko and Pastur, 1967). Readers may refer to Jin et al. (2014) and Wang et al. (2015) for more details about the model.

To estimate $k$ and $q$, the following method can be employed. First, note that when $\tau=0$ and $\textrm{Cov}(\bbf_t)=\Sigma_\bbf$, the population covariance matrix of $\bbR_{t}$ is a \emph{spiked population model} (Johnstone (2001), Baik and Silverstein (2006), Bai and Yao (2008)) with $k(q+1)$ spikes. Therefore, $k(q+1)$ can be estimated by counting the number of eigenvalues of $\Phi_n(0)$ that are larger than some phase transition point. Next, the separated estimation of $k$ and $q$ can be achieved by investigating the spectral property of $\bbM_{n}(\tau)$ for general $\tau \ge 1$, using the fact that the number of eigenvalues of $\Phi_{n}(\tau)$ that lie outside the support of the LSD of $\bbM_n(\tau)$ at lags $1 \leq \tau \leq q $ is different from that at lags $\tau>q$. Thus, the estimates of $k$ and $q$ can be separated by counting the number of eigenvalues of $\Phi_{n}(\tau)$ that lie outside the support of the LSD of $\bbM_n(\tau)$ from $\tau=0,1,2,\cdots,q,q+1,\cdots.$

Note that for the above method to work, the LSD of $\bbM_n(\tau)$ for general $\tau \ge 1$ must be known. This is derived in Jin et al. (2014).
Moreover, it is required that no eigenvalues outside the the support of the LSD of $\bbM_n(\tau)$ so that if an eigenvalue of $\Phi_{n}(\tau)$ goes out of the support of the LSD of $\bbM_n(\tau)$, it must come from the signal part. Wang et al. (2015) proved such phenomenon theoretically. Both results are included in Section 2 for readers' reference.

The rest of the paper is structured as follows: Some known results are given in Section 2. Section 3 presents truncation of variables and Section 4 estimates $k(q+1)$. The estimation of $q$ is provided in Section 5, from the which the estimation of $k$ can also be obtained. Section 6 discusses the case when the variance of the noise part is unknown. A simulation study is shown in Section 7 and some proofs are presented in Appendix.

Regrading the norm used in this paper, the norm applied to a vector is the usual Euclidean norm, with notation $\|*\|$. For a matrix, two kinds of norm have been used. The operator norm, denoted by $\|*\|_o$, is the largest singular value. For matrices of fixed dimension, the Kolmogorov norm, defined as the largest absolute value of all the entries, has been used, with notation $\|*\|_K$.

\section{Some known results}
\setcounter{equation}{0}
\def\theequation{\thesection.\arabic{equation}}
\setcounter{subsection}{0}
In this section, we present some known results.
\begin{lemma}\label{Burkholder}(Burkholder (1973)). Let $\{X_k\}$ be a complex martingale difference sequence with respect to the increasing $\sigma$-fields $\{\mathcal{F}_n\}$. Then, for $p\ge2$, we have
\bqn
\emph{\rE}|\sum X_k|^p\le K_p\Big(\emph{\rE}\big(\sum\emph{\rE}(|X_k|^2|\mathcal{F}_{k-1})\big)^{p/2}+\emph{\rE}\sum|X_k|^p\Big).
\eqn
\end{lemma}
\begin{lemma}\label{A1} (Lemma A.1 of Bai and Silverstein (1998)). For $X=(X_1,\cdots,X_n)'$ i.i.d. standardized (complex) entries,
$\bbB$ $n\times n$ Hermitian nonnegative definite matrix, we have, for any $p\ge1$,
\bqn
\emph{\rE}|\bbX^*\bbB\bbX|^p\le K_p\Big(\big(\emph{\rtr}\bbB\big)^p+rE|X_1|^{2p}\emph{\rtr}\bbB^p\Big),
\eqn
where $K_p$ is a constant depending on $p$ only.
\end{lemma}
\begin{lemma} \label{Jin}(Jin et al. (2014)). Assume:
 \begin{description}
   \item(a) $\tau \ge1$ is a fixed integer.
   \item(b) $\bbe_{k}=(\ep_{1k},\cdots,\ep_{nk})'$, $k=1,2,...,T+\tau$, are $n$-dimensional vectors of independent standard complex components with $\sup_{1\le i\le n,1\le t\le T+\tau}\rE|\ep_{it}|^{2+\delta}\le M <\infty$ for some $\delta\in(0,2]$, and for any $\eta>0$, \bqa \label{eta}  \frac{1}{\eta^{2+\delta}n T}\sum_{i=1}^n\sum_{t=1}^{T+\tau}\rE(|\ep_{it}|^{2+\delta}I(|\ep_{it}|\geq \eta T^{1/(2+\delta)}))=o(1).
\eqa
   \item(c) $n/(T+\tau)\rightarrow c>0$ as $n,T\to \infty$.
   \item(d) $\bbM_n(\tau)=\sum_{k=1}^{T}(\bgma_{k}\bgma_{k+\tau}^{*}+\bgma_{k+\tau}\bgma_{k}^{*}),$ where $\bgma_{k}=\frac{1}{\sqrt {2T}}\bbe_{k}$.
 \end{description}
 Then as $n,T\rightarrow\infty$, $F^{\bbM_n(\tau)}\overset {D}{\rightarrow} F_c$ a.s. and $F_c$ has a density function given by
 \begin{center}
 $\phi_c(x)=\frac{1}{2c\pi}\sqrt{\frac{y_0^2}{1+y_0}-(\frac{1-c}{|x|}+\frac{1}{\sqrt{1+y_0}})^2}$, $ |x|\leq a,$
 \end{center}
 where
\begin{eqnarray*}
a=\left\{\begin{array}{cc}\frac{(1-c)\sqrt{1+y_1}}{y_1-1},& c \neq 1,\\
2,& c=1,
\end{array}\right.
\end{eqnarray*}
$y_0$ is the largest real root of the equation: $y^3-\frac{(1-c)^2-x^2}{x^2}y^2-\frac{4}{x^2}y-\frac{4}{x^2}=0$ and $y_1$ is the only real root of the equation: \bqa\label{ey1}((1-c)^2-1)y^3+y^2+y-1=0\eqa such that $y_1>1$ if $c<1$ and $y_1\in (0,1)$ if $c>1$.
Further, if $c>1$, then $F_c$ has a point mass $1-1/c$ at the origin. Note that as long as $\tau\ge1$, $F_c$ does not depend on $\tau.$
\end{lemma}
\begin{lemma}(Bai and Wang (2015)). \label{Bai} Theorem \ref{Jin} still holds with the $2+\delta$ moment condition weakened to 2nd moment.
\end{lemma}
\begin{lemma} \label{Wang} (Wang et al. (2015)). Assume:
 \begin{description}
   \item(a) $\tau \ge1$ is a fixed integer.
   \item(b) $\bbe_{k}=(\ep_{1k},\cdots,\ep_{nk})'$, $k=1,2,...,T+\tau$, are $n$-vectors of independent standard complex components with $\sup_{i, t}\rE|\ep_{it}|^4\le M$ for some $M>0$.
   \item(c) There exist $K>0$ and a random variable $X$ with finite fourth order moment such that, for any $x>0$, for all $n,T$
       \bqa \label{eta}
       \frac{1}{nT}\sum_{i=1}^n\sum_{t=1}^{T+\tau}\rP(|\ep_{it}|> x)\le K\rP(|X|>x).
       \eqa
   \item(d) $c_n\equiv n/T\rightarrow c>0$ as $n\to \infty$.
   \item(e) $\bbM_n=\sum_{k=1}^{T}(\bgma_{k}\bgma_{k+\tau}^{*}+\bgma_{k+\tau}\bgma_{k}^{*}),$ where $\bgma_{k}=\frac{1}{\sqrt {2T}}\bbe_{k}$.
   \item(f) The interval [a,b] lies outside the support of $F_{c}$, where $F_c$ is defined as in Lemma \ref{Jin}.
 \end{description}
 Then $\rP\left(\mbox{no eigenvalues of $\bbM_n$ appear in $[a,b]$ for all large n}\right)=1$.
\end{lemma}
\section{\large Truncation, centralization and standardization of variables}
As proved in Wang et al.(2015), we may assume that the $\ep_{ij}$'s satisfy the conditions that
\bqa\label{trun}
|\ep_{ij}|\le C,~\rE \ep_{ij}= 0,
~\rE |\ep_{ij}|^2=1,~\rE |\ep_{ij}|^4<M
\eqa
for some $C, M>0$.

For the truncation of variables in $\bbF^\tau$, first note that for a random variable $X$ with $\rE|X|^4<\infty$, we have
$\sum_{\ell=1}^\infty2^\ell\rP\left(|X|>2^{\ell/4}\right)<\infty$. 
Given $\rE|\bbF^{\tau}_{ij}|^4<\infty, i=1,\cdots, k(q+1), j=1,\cdots, T$, define
$
\hat\bbF^{\tau}_{ij}=\left\{\begin{array}{cc}\bbF^{\tau}_{ij},& |\bbF^{\tau}_{ij}|<T^{1/4},\\
0,& \mbox{otherwise},
\end{array}\right.
$
$\hat\bbF^{\tau}=\big(\hat\bbF^{\tau}_{ij}\big)$ and $$\hat\bbPhi_{n}(\tau)=\frac1{2T}[(\bgL\hat\bbF^{\tau}+\boldsymbol{\bbe}^{\tau})(\bgL\hat\bbF^{0}+\boldsymbol{\bbe}^{0})^*
+(\bgL\hat\bbF^{0}+\boldsymbol{\bbe}^{0})(\bgL\hat\bbF^{\tau}+\boldsymbol{\bbe}^{\tau})^*].$$
Then we have
\bqn
&&\rP\left(\bbPhi_{n}(\tau)\ne\hat\bbPhi_{n}(\tau), i.o.\right)\\
&=&\rP\left(\bbF^{\tau}\ne\hat\bbF^{\tau}, i.o.\right)\\
&=&\rP\left(\bigcap_{L=1}^\infty\bigcup_{T=L}^\infty\bigcup_{\begin{subarray}{c}i\le k(q+1)\\j\le T\end{subarray}}\{|\bbF^{\tau}_{ij}|\ge T^{1/4}\}\right)\\
&\le&\lim_{L\to\infty}\sum_{\ell=L}^{\infty}\rP\left(\bigcup_{T=2^\ell+1}^{2^{\ell+1}}\bigcup_{\begin{subarray}{c}i\le k(q+1)\\j\le 2^{\ell+1}\end{subarray}}\{|\bbF^{\tau}_{ij}|\ge2^{\ell/4}\}\right)\\
&\le&\lim_{L\to\infty}\sum_{\ell=L}^{\infty}\rP\left(\bigcup_{\begin{subarray}{c}i\le k(q+1)\\j\le 2^{\ell+1}\end{subarray}}\{|\bbF^{\tau}_{ij}|\ge2^{\ell/4}\}\right)\\
&\le&k(q+1)\lim_{L\to\infty}\sum_{\ell=L}^{\infty}2^{\ell+1}\rP\left(|\bbF^{\tau}_{11}|\ge2^{\ell/4}\right)\\
&\to&0.
\eqn
This completes the proof of truncation. Centralization and standardization can be justified in the same way as in Appendix A of Wang et al. (2015). In what follows, we may assume that
\bqn
|\bbF^{\tau}_{ij}|< T^{1/4},~\rE\bbF^{\tau}_{ij}= 0,~\rE|\bbF^{\tau}_{ij}|^2=1, ~\rE |\bbF^{\tau}_{ij}|^4<M
\eqn
for some $M>0$.
\section{Estimation of $k(q+1)$}
\setcounter{equation}{0}
\def\theequation{\thesection.\arabic{equation}}
\setcounter{subsection}{0}
In this section, we will estimate $k(q+1)$ by an investigation of the limiting properties of eigenvalues of $\bbPhi_n(0)$. For simplicity, rewrite $\bbPhi_n(0)=\bbPhi(0), \bbF^0=\bbF$ and $\bbe^0=\bbe$. With these notations, we have $\bbPhi(0)=\frac1T(\bgL\bbF+\bbe)(\bgL\bbF+\bbe)^*$ and $\bbM(0)=\frac1T\bbe\bbe^*$.
When $\bgL=\textbf{0}$, $\bbPhi(0)$ reduces to $\bbM(0)$, which is a standard sample covariance matrix and thus its ESD tends to the famous MP law (Mar\v{c}enko and Pastur, 1967).

Suppose $\ell$ is an eigenvalue of $\bbPhi(0)$, then we have
\bqa
0=\det\left|\ell\bbI-\bbPhi(0)\right|=\det\left|\ell\bbI-\bbM(0)-\frac1T \bgL\bbF\bbe^*-\frac1T\bbe\bbF^*\bgL^*-\frac1T\bgL\bbF\bbF^*\bgL^*\right|.
\label{eqn1}
\eqa
Let $\bbB=(\bbB_1:\bbB_2)$ be an $n\times n$ orthogonal matrix such that $\bbB_1=\bgL(\bgL^*\bgL)^{-1/2}$ and thus $\bgL^*\bbB_2=\textbf{0}_{k(q+1)\times(n-k(q+1))}$. Then (\ref{eqn1}) is equivalent to
\bqa
\det\begin{vmatrix} \ell\bbI_{k(q+1)}-\frac1T\bbB_1^*(\bgL\bbF+\bbe)(\bbF^*\bgL^*+\bbe^*)\bbB_1&-\frac1T\bbB_1^*(\bgL\bbF+\bbe)\bbe^*\bbB_2\cr
-\frac1T\bbB_2^*\bbe(\bbF^*\bgL^*+\bbe^*)\bbB_1&\ell\bbI_{n-k(q+1)}-\frac1T\bbB_2^*\bbe\bbe^*\bbB_2\cr\end{vmatrix}=0
\eqa
If we further assume that $\ell$ is not an eigenvalue of $\frac1T\bbB_2^*\bbe\bbe^*\bbB_2$, then we have
\bqa
\det|\bbI_{k(q+1)}-\frac1T\bbB_1^*(\bgL\bbF+\bbe)\bbD^{-1}(\ell)(\bbF^*\bgL^*+\bbe^*)\bbB_1|=0,
\label{eq3}
\eqa
where $\bbD(\ell)=\ell\bbI_T-\frac1T\bbe^*\bbB_2\bbB_2^*\bbe$. Denote $\bbH(\ell)=\ell\bbI_T-\frac1T\bbe^*\bbe$, then we obtain
\bqa\label{eq6}
\frac1T\bbB_1^*\bbe\bbD^{-1}(\ell)\bbe^*\bbB_1
=\Big(\bbI+\frac1T\bbB_1^*\bbe\bbH^{-1}(\ell)\bbe^*\bbB_1\Big)^{-1}\frac1T\bbB_1^*\bbe\bbH^{-1}(\ell)\bbe^*\bbB_1.
\eqa
Next, we have
\bqa\label{eq7}
\frac1T\bbB_1^*\bbe\bbH^{-1}(\ell)\bbe^*\bbB_1&=&\frac1T\bbB_1^*\Big(-T\bbI+\ell T(\ell\bbI_n-\bbM(0))^{-1}\Big)\bbB_1\non
&=&-\bbI_{k(q+1)}+\ell\bbB_1^*(\ell\bbI_n-\bbM(0))^{-1}\bbB_1.
\eqa
Substitute (\ref{eq7}) back to (\ref{eq6}), and we have
\bqn
&&\frac1T\bbB_1^*\bbe\bbD^{-1}(\ell)\bbe^*\bbB_1\\
&=&\bbI_{k(q+1)}-\Big(\ell\bbB_1^*(\ell\bbI_n-\bbM(0))^{-1}\bbB_1\Big)^{-1}=\bbI_{k(q+1)}+\frac1\ell\Big(\bbB_1^*(\bbM(0)-\ell\bbI_n)^{-1}\bbB_1\Big)^{-1}
\eqn
Write $\bbB_1=(\bbb_1,\cdots,\bbb_{k(q+1)})$, then we have $\|\bbb_i\|=1$. By Lemma 6 in Bai, Liu and Wong (2011), we have
\bqn
\bbb_i^*(\bbM(0)-\ell\bbI_n)^{-1}\bbb_i\to m,  \quad a.s.
\eqn
and for $i\ne j$
\bqn
\bbb_i^*(\bbM(0)-\ell\bbI_n)^{-1}\bbb_j\to 0, \quad a.s.,
\eqn
where $m=m(\ell)=\lim_{n\to\infty}\frac1n\rtr(\bbM(0)-\ell\bbI_n)^{-1}$ is the Stieltjes transform of the sample covariance with ratio index $c=\lim_{n\to\infty}\frac nT$.

By Lemma 3.11 of Bai and Silverstein (2010), we have $m$ satisfying
\bqa\label{ml}
m(\ell)=\frac{1-c-\ell+\sqrt{(1-\ell-c)^2-4\ell c}}{2c\ell}.
\eqa
Therefore, we obtain
\bqn
\frac1T\bbB_1^*\bbe\bbD^{-1}(\ell)\bbe^*\bbB_1\to\big(1+\frac1{\ell m}\big)\bbI_{k(q+1)}.
\eqn
Next, we want to show, with probability 1 that $$\frac1T\bbB_1^*\bgL\bbF\bbD^{-1}(\ell)\bbe^*\bbB_1\to \textbf{0}$$ and $$\frac1T\bbB_1^*\bbe\bbD^{-1}(\ell)\bbF^*\bgL^*\bbB_1\to \textbf{0}.$$
Note that
\bqn
&&\frac1T\bbB_1^*\bgL\bbF\bbD^{-1}(\ell)\bbe^*\bbB_1\\
&=&\frac1T\big(\bbI+\frac1T\bbB_1^*\bbe\bbH^{-1}(\ell)\bbe^*\bbB_1\big)^{-1}\bbB_1^*\bgL\bbF\bbH^{-1}(\ell)\bbe^*\bbB_1\\
&=&\frac1{\ell T}\Big(\bbB_1^*(\bbM(0)-\ell\bbI_N)^{-1}\bbB_1\Big)^{-1}\bbB_1^*\bgL\bbF\bbH^{-1}(\ell)\bbe^*\bbB_1\\
&=&\frac1{\ell T}\Big(\bbB_1^*(\bbM(0)-\ell\bbI_N)^{-1}\bbB_1\Big)^{-1}\bbB_1^*\bgL\bbF\big(\ell\bbI_T-\frac1T\bbe^*\bbe\big)^{-1}\bbe^*\bbB_1.
\eqn
Recall $\bbM(0)=\frac1T\bbe\bbe^*$. Fix $\delta>0$ and let event $\mathcal{A}=\{\lambda_{\max}(\bbM(0))\le(1+\sqrt c)^2+\delta\}$ and $\mathcal{A}^c$ be the complement. By Theorem 5.9 of Bai and Silverstein (2010), we have $\rP(\mathcal{A}^c)=o(n^{-t})$ for any $t>0$.

Suppose $\ell$ is an eigenvalue of $\bbPhi(0)$ larger than $(1+\sqrt c)^2+2\delta$. By the fact that $\bbM(0)$ and $\frac1T\bbe^*\bbe$ have the same set of nonzero eigenvalues, we have, under $\mathcal{A}$, that $\|\ell\bbI_T-\frac1T\bbe^*\bbe\|_o\ge \ell-\|\frac1T\bbe^*\bbe\|_o\ge\delta>0$, and hence $\|\big(\ell\bbI_T-\frac1T\bbe^*\bbe\big)^{-1}\|_o\le\frac 1\delta$.

Therefore, for any $\ep>0$, we have
\bqn
&&\rP(\|\frac1T\bbF\big(\ell\bbI_T-\frac1T\bbe^*\bbe\big)^{-1}\bbe^*\bbB_1\|_K\ge\ep)\\
&=&\rE\Big(\rP(\|\frac1T\bbF\big(\ell\bbI_T-\frac1T\bbe^*\bbe\big)^{-1}\bbe^*\bbB_1\|_K\ge\ep)\Big|\bbe\Big)\\
&\le&\rE\Big(\rP(\|\frac1T\bbF\big(\ell\bbI_T-\frac1T\bbe^*\bbe\big)^{-1}\bbe^*\bbB_1\|_K\ge\ep)\Big|\bbe,\mathcal{A}\Big)+\rP(\mathcal{A}^c).
\eqn
Write $\bbF=(\tilde\bbF_1,\cdots,\tilde\bbF_{k(q+1)})'$.
For the first term, by Lemma \ref{A1}, we have
\bqn
&&\rE\Big(\rP(\|\frac1T\bbF\big(\ell\bbI_T-\frac1T\bbe^*\bbe\big)^{-1}\bbe^*\bbB_1\|_K\ge\ep)\Big|\bbe,\mathcal{A}\Big)\\
&\le&\frac1{\ep^{4r}T^{4r}}\rE\rE\Big(\|\bbF\big(\ell\bbI_T-\frac1T\bbe^*\bbe\big)^{-1}\bbe^*\bbB_1\|_K^{4r}\Big|\bbe,\mathcal{A}\Big)\\
&\le&\frac1{\ep^{4r} T^{4r}}\rE\rE\Big[\Big(\rtr\bbF\big(\ell\bbI_T-\frac1T\bbe^*\bbe\big)^{-1}\bbe^*\bbB_1\bbB_1^*\bbe\big(\ell\bbI_T-\frac1T\bbe^*\bbe\big)^{-1}\bbF^*\Big)^{2r}\Big|\bbe,\mathcal{A}\Big]\\
&=&\frac1{\ep^{4r}T^{4r}}\rE\rE\Big[\Big(
\sum_{i=1}^{k(q+1)}\tilde\bbF_i^*\big(\ell\bbI_T-\frac1T\bbe^*\bbe\big)^{-1}\bbe^*\bbB_1\bbB_1^*\bbe\big(\ell\bbI_T-\frac1T\bbe^*\bbe\big)^{-1}\tilde\bbF_i\Big)^{2r}\Big|\bbe,\mathcal{A}\Big]\\
&\le&\frac {[k(q+1)]^{2r-1}}{\ep^{4r}T^{4r}}\rE\Big[
\sum_{i=1}^{k(q+1)}\rE\Big(\tilde\bbF_i^*\big(\ell\bbI_T-\frac1T\bbe^*\bbe\big)^{-1}\bbe^*\bbB_1\bbB_1^*\bbe\big(\ell\bbI_T-\frac1T\bbe^*\bbe\big)^{-1}\tilde\bbF_i\Big)^{2r}\Big|\bbe,\mathcal{A}\Big]\\
&\le&\frac {K_{2r}[k(q+1)]^{2r-1}}{\ep^{4r}T^{4r}}\rE\Big[\sum_{i=1}^{k(q+1)}\Big(
\big[\rtr\big(\ell\bbI_T-\frac1T\bbe^*\bbe\big)^{-1}\bbe^*\bbB_1\bbB_1^*\bbe\big(\ell\bbI_T-\frac1T\bbe^*\bbe\big)^{-1}
\big]^{2r}+\\
&&\rE|\bbF_{11}|^{4r}\rtr[(\ell\bbI_T-\frac1T\bbe^*\bbe\big)^{-1}\bbe^*\bbB_1\bbB_1^*\bbe\big(\ell\bbI_T-\frac1T\bbe^*\bbe\big)^{-1}]^{2r}\Big)\Big| \mathcal{A}\Big]\\
&\le&\frac{K_{2r}[k(q+1)]^{2r}}{\delta^{4r}\ep^{4r}T^{4r}}\rE\Big[
\big(\rtr\bbe^*\bbB_1\bbB_1^*\bbe\big)^{2r}+\rE|\bbF_{11}|^{4r}\rtr(\bbe^*\bbB_1\bbB_1^*\bbe)^{2r}\Big| \mathcal{A}\Big]\\
&\le&\frac {K_{2r}[k(q+1)]^{2r}[(1+\sqrt c)^2+\delta]^{2r}}{\delta^{4r}\ep^{4r}T^{2r}}\rE\Big((\rtr\bbB_1^*\bbB_1)^{2r}+\rE|\bbF_{11}|^{4r}\rtr(\bbB_1^*\bbB_1)^{2r}\Big)\\
&\le&\frac {K_{2r}[k(q+1)]^{2r}[(1+\sqrt c)^2+\delta]^{2r}}{\delta^{4r}\ep^{4r}T^{2r}}\big([k(q+1)]^{2r}+\rE|\bbF_{11}|^4T^{\frac{4r-4}4}k(q+1)\big)
\eqn
which is summable for $r\ge1$.\\
Hence, we have shown with probability 1 that $$\frac1T\bbB_1^*\bgL\bbF\bbD^{-1}(\ell)\bbe^*\bbB_1\to \textbf{0}.$$
Similarly, we have with probability 1 that $$\frac1T\bbB_1^*\bbe\bbD^{-1}(\ell)\bbF^*\bgL^*\bbB_1\to \textbf{0}.$$
Therefore, substituting into (\ref{eq3}), we have
\bqa
\det|\frac1T\bbB_1^*\bgL\bbF\bbD^{-1}(\ell)\bbF^*\bgL^*\bbB_1+\frac1{\ell m(\ell)}\bbI_{k(q+1)}|\to 0.
\label{eq4}
\eqa
Using Bai, Liu and Wong (2011) again, we have the diagonal elements
of the matrix $T^{-1}\bbF\bbD^{-1}(\ell)\bbF^*$ tend to $-\underline{m}(\ell)$ and the off diagonal elements tend to 0. Here $\underline{m}(\ell)$ is the Stieltjies transform of the LSD of $\frac1T\bbe^*\bbe$ and satisfies
$$\underline{m}(\ell)=-\frac{1-c}\ell+cm(\ell).$$
Thus, if $\bgL^*\bgL\to \bbQ$, then (\ref{eq4}) can be further simplified as
\bqa\label{eq5}
\det|-\bbQ\underline{m}(\ell)+\frac1{\ell m(\ell)}\bbI_{k(q+1)}|=0.
\eqa
If $\ga$ is an eigenvalue of $\bbQ$, and there is an $\ell$ belonging to the complement of the support of the LSD of $\bbM(0)$ such that
$\ga=\frac1{\ell m(\ell)\underline{m}(\ell)}$, then $\ell$ is a solution of (\ref{eq5}).

From (\ref{ml}), we have
\bqn
c\ell m^2(\ell)-(1-c-\ell)m(\ell)+1=0,
\eqn
which implies
\bqn
\ell m(\ell)\underline{m}(\ell)&=&\ell m(\ell)\big(-\frac{1-c}\ell+cm(\ell)\big)\\
&=&-(1-c)m(\ell)+c\ell m^2(\ell)\\
&=&-(1-c)m(\ell)+(1-c-\ell)m(\ell)-1\\
&=&-\ell m(\ell)-1\\
&=&-\frac{1-c-\ell+\sqrt{(1-\ell-c)^2-4\ell c}}{2c}-1\\
&=:&g(\ell).
\eqn
It is easy to verify that $g'(\ell)<0$, implying that $\ell m(\ell)\underline{m}(\ell)$ is decreasing. Also note that $\ell m(\ell)\underline{m}(\ell)=\frac1{\sqrt c}$ when
$\ell=(1+\sqrt{c})^2$. Therefore, if $\ga=\frac1{\ell m(\ell)\underline{m}(\ell)}>\sqrt c$, then we have $\ell>(1+\sqrt{c})^2$. This recovers the result of Baik and Silverstein (2006). Note that $[(1-\sqrt c)^2,(1+\sqrt c)^2]$ is the support of the MP law. Hence, if all the eigenvalues of $\bbQ$ are greater than $\sqrt c$, we have $k(q+1)$ sample eigenvalues of $\bbPhi_n(0)$ goes outside the right boundary of the support of the MP law. Note that although the distribution of rest $n-k(q+1)$ sample eigenvalues follows the MP law with the largest sample eigenvalue converging to the right boundary, there is still a positive probability that the largest sample eigenvalue goes beyond the right boundary. Therefore, to completely separate the $k(q+1)$ spiked sample eigenvalues from the rest, the threshold is set as $(1+\sqrt c)^2(1+2n^{-2/3})$.  In other words, $k(q+1)$ can be estimated by the number of sample eigenvalues of $\bbPhi_n(0)$ greater than $(1+\sqrt c)^2(1+2n^{-2/3})$.
\begin{remark}
For factor models, the loading matrix is unknown. This, however, is not a concern in our estimation because compared with the noise matrix, the loading matrix is denominating, making the condition easily satisfied that all the eigenvalues of $\bbQ$ are greater than $\sqrt c$.
\end{remark}
\begin{remark}
The rationale of choosing $(1+2n^{-2/3})$ as the buffering factor of the criterion is that, according to Tracy-Widom law, the quantity of a non-spiked eigenvalue larger than $(1+\sqrt c)^2$ has an order of $n^{-2/3}$. Therefore, it is good enough for us to choose $(1+2n^{-2/3})$ to completely separate the spikes and the bulk eigenvalues.
\end{remark}
\section{Estimation of $q$}
\setcounter{equation}{0}
\def\theequation{\thesection.\arabic{equation}}
\setcounter{subsection}{0}
Next, we want to split $k$ and $q$. Let $\tau\ge1$ be given and assume that $\ell$ is an eigenvalue of $\bbPhi_n(\tau)$. For simplicity, write $\bbM_n(\tau)=\bbM$ and for $t=1,2,\cdots,T$, define $\bbF_t=(\bbf_t,\bbf_{t-1},\cdots,\bbf_{t-q})'$ such that $\bbR_{t}=\bgL\bbF_t+\bbe_{t}$.
Then we have
\bqa
0&=&\det\left|\ell\bbI-\bbPhi_n(\tau)\right|\non
&=&\det\Big|\ell\bbI-\bbM-\frac{1}{2T}
\sum_{j=1}^{T}\Big(\bgL \bbF_j \bbF_{j+\tau}^*\bgL^*+\bgL \bbF_{j+\tau}^*\bbF_j^*\bgL^*\non
&&+\bbe_j \bbF_{j+\tau}^*\bgL^*+\bgL \bbF_{j+\tau}^*\bbe_j^*+\bbe_{j+\tau}^*\bbF_j^*\bgL^*+\bgL \bbF_j \bbe_{j+\tau}^*\Big)\Big|.
\label{eq0}
\eqa
Define $\bbB, \bbB_1$ and $\bbB_2$ the same as in the last section. Multiplying $\bbB^*$ from left and $\bbB$ from right to the above matrix and by $\bgL^*\bbB_2=\textbf{0}_{k(q+1)\times(n-k(q+1))}$, we have (\ref{eq0}) equivalent to
\bqn
0=\det\begin{vmatrix}\ell\bbI_{k(q+1)}-\bbS_{11}&-\bbS_{12}\cr-\bbS_{21}&\ell\bbI_{n-k(q+1)}-\bbS_{22}\end{vmatrix}=\det|\ell\bbI-\bbS_{22}|\det|\ell\bbI-\bbK_n(\ell)|,
\eqn
where
\bqn
\bbS_{11}&=&\frac1{2T}\sum_{j=1}^T\bbB_1^*[(\bgL \bbF_j+\bbe_j)(\bgL \bbF_{j+\tau}+\bbe_{j+\tau})^*+(\bgL \bbF_{j+\tau}+\bbe_{j+\tau})(\bgL \bbF_j+\bbe_j)^*]\bbB_1\\
\bbS_{12}&=&\frac1{2T}\sum_{j=1}^T\bbB_1^*(\bgL \bbF_j\bbe_{j+\tau}^*+\bgL \bbF_{j+\tau}\bbe_j^*)\bbB_2+\bbB_1^*\bbM\bbB_2\\
\bbS_{21}&=&\bbS_{12}^*\\
\bbS_{22}&=&\bbB_2^*\bbM\bbB_2\\
\bbK_n(\ell)&=&\bbS_{11}+\bbS_{12}(\ell\bbI_{n-k(q+1)}-\bbS_{22})^{-1}\bbS_{21}.
\eqn
Therefore, if $\ell$ is not an eigenvalue of $\bbS_{22}$, by the factorization above, $\ell$ must be an eigenvalue of $\bbK_n(\ell)$, i.e.
$\det|\bbK_n(\ell)-\ell\bbI|=0$.

Denote $\bbW=\frac1{2T}\sum_{j=1}^T(\bbF_j\bbe_{j+\tau}^*+\bbF_{j+\tau}\bbe_j^*)$.
By the assumptions of $\bbe_t$'s and $\bbF_t$'s, the random vector $\{\bbF_j\bbe_{j+\tau}^*+\bbF_{j+\tau}\bbe_j^*, j\ge1\}$ is $(q+1)$-dependent (see Page 224, Chung 2001). It then follows with probability 1 that
\bqn
\bbW\bbB_1&=&o(\textbf{1})\\
\bbB_1^*\bbW^*&=&o(\textbf{1})\\
\frac1{2T}\sum_{j=1}^T(\bbF_j\bbF_{j+\tau}^*+\bbF_{j+\tau}\bbF_j^*)&=&\bbH(\tau)+o(\textbf{1}),
\eqn
where
\bqn
\bbH(\tau)&=&\begin{pmatrix}
  0 & \cdots & 1 & \cdots & 0 \\
  \vdots & \ddots & 0 &  1 &  \vdots\\
 1 &  0 &\ddots&0 & 1  \\
  \vdots & 1& 0 & \ddots & \vdots   \\
  0 &\cdots & 1& \cdots & 0  \\
 \end{pmatrix},
\eqn
is of dimension $k(q+1)\times k(q+1)$ with two bands of 1's of $k\tau$-distance from the main diagonal.
Therefore, we have a.s.
\bqn
\bbS_{11}&=&\bbB_1^*\bgL\bbH(\tau)\bgL^*\bbB_1+\bbB_1^*\bbM\bbB_1+o(\textbf{1})\\
\bbS_{12}&=&\bbB_1^*(\bbM+\bgL\bbW)\bbB_2\\
\bbS_{21}&=&\bbB_2^*(\bbM+\bbW^*\bgL^*)\bbB_1.
\eqn
Subsequently, we have a.s.
\bqn
\bbK_n(\ell)=\bbB_1^*\bgL\bbH(\tau)\bgL^*\bbB_1+\bbB_1^*\bbM\bbB_1
+\bbB_1^*(\bbM+\bgL\bbW)\bbB_2\Big(\ell\bbI-\bbB_2^*\bbM\bbB_2\Big)^{-1}\bbB_2^*(\bbM+\bbW^*\bgL^*)\bbB_1+o(\textbf{1}).
\eqn
Note that
\bqn
&&\bbB_1^*\bbM\bbB_1+\bbB_1^*\bbM\bbB_2\Big(\ell\bbI-\bbB_2^*\bbM\bbB_2\Big)^{-1}\bbB_2^*\bbM\bbB_1\\
&=&\bbB_1^*\bbM\bbB_1+\bbB_1^*\bbM\frac1\ell\bbB_2\bbB_2^*\bbM(\bbI-\frac1\ell\bbB_2\bbB_2^*\bbM)^{-1}\bbB_1\\
&=&\bbB_1^*\bbM(\bbI-\frac1\ell\bbB_2\bbB_2^*\bbM)^{-1}\bbB_1\\
&=&\ell\bbB_1^*\bbM(\ell\bbI-\bbM+\bbB_1\bbB_1^*\bbM)^{-1}\bbB_1\\
&=&\ell\bbI-\ell\Big(\bbI+\bbB_1^*\bbM\big(\ell\bbI-\bbM\big)^{-1}\bbB_1\Big)^{-1}\\
&=&\ell\bbI-\Big(\bbB_1^*\big(\ell\bbI-\bbM\big)^{-1}\bbB_1\Big)^{-1}
\eqn
and
\bqn
&&\bbW\bbB_2\Big(\ell\bbI-\bbB_2^*\bbM\bbB_2\Big)^{-1}\bbB_2^*\bbW^*\\
&=&\frac1{2T}\sum_{i,j=1}^{T}(\bbF_i\bgma_{i+\tau}^*+\bbF_{i+\tau}\bgma_i^*)\bbB_2\Big(\ell\bbI-\bbB_2^*\bbM\bbB_2\Big)^{-1}\bbB_2^*(\bgma_{j+\tau}\bbF_j^*+\bgma_j\bbF_{j+\tau}^*)\\
&=&\frac1{2T}\sum_{j=1}^{T}\Big[\bbF_j\bgma_{j+\tau}^*\bbB_2\Big(\ell\bbI-\bbB_2^*\bbM\bbB_2\Big)^{-1}\bbB_2^*\bgma_{j+\tau}\bbF_j^*
+\bbF_{j+\tau}\bgma_j^*\bbB_2\Big(\ell\bbI-\bbB_2^*\bbM\bbB_2\Big)^{-1}\bbB_2^*\bgma_j\bbF_{j+\tau}^*\Big]+\\
&&\frac1{2T}\sum_{j=1}^{T}\Big[\bbF_j\bgma_{j+\tau}^*\bbB_2\Big(\ell\bbI-\bbB_2^*\bbM\bbB_2\Big)^{-1}\bbB_2^*\bgma_j\bbF_{j+\tau}^*
+\bbF_{j+\tau}\bgma_j^*\bbB_2\Big(\ell\bbI-\bbB_2^*\bbM\bbB_2\Big)^{-1}\bbB_2^*\bgma_{j+\tau}\bbF_j^*\Big]+\\
&&\frac1{2T}\sum_{j=1}^{T}(\bbF_{j+\tau}\bgma_{j+2\tau}^*+\bbF_{j+2\tau}\bgma_{j+\tau}^*)\bbB_2\Big(\ell\bbI-\bbB_2^*\bbM\bbB_2\Big)^{-1}\bbB_2^*(\bgma_{j+\tau}\bbF_j^*+\bgma_j\bbF_{j+\tau}^*)+\\
&&\frac1{2T}\sum_{j=1}^{T}(\bbF_{j-\tau}\bgma_{j}^*+\bbF_{j}\bgma_{j-\tau}^*)\bbB_2\Big(\ell\bbI-\bbB_2^*\bbM\bbB_2\Big)^{-1}\bbB_2^*(\bgma_{j+\tau}\bbF_j^*+\bgma_j\bbF_{j+\tau}^*)+\\
&&\frac1{2T}\sum_{i,j=1 \atop i\ne j,i\ne j\pm\tau}^{T}(\bbF_i\bgma_{i+\tau}^*+\bbF_{i+\tau}\bgma_i^*)\bbB_2\Big(\ell\bbI-\bbB_2^*\bbM\bbB_2\Big)^{-1}\bbB_2^*(\bgma_{j+\tau}\bbF_j^*+\bgma_j\bbF_{j+\tau}^*)\\
&=:&\bbP_1+\bbP_2+\bbP_3+\bbP_4+\bbP_5.
\eqn
Next, we give a lemma on the quadratic form of $\bgma_{j}$.
\begin{lemma}\label{Cp}Let $i,j\in \mathbb{N}$ be given, we have almost surely and uniformly in $i$ and $j$ that
\bqn
\bgma_i^*\bbB_2\Big(\ell\bbI-\bbB_2^*\bbM\bbB_2\Big)^{-1}\bbB_2^*\bgma_j\to
\left\{\begin{array}{l}\frac{-\frac{cm}2}{1-\frac{c^2m^2}{2x_1}}\Big(\frac{-\frac{cm}2}{1-\frac{c^2m^2}{4x_1}}\Big)^p\equiv C_p, \quad i=j\pm p\tau\\
0, \quad \mbox{otherwise}.
\end{array}\right.
\eqn
\end{lemma}
The proof of the lemma is postponed in the Appendix.
First, we have
\bqn
&&\rE(\bbP_1)\\
&=&\rE[\rE(\bbP_1|\bgma_1,\cdots,\bgma_{T+\tau})]\\
&=&\frac1{2T}\sum_{j=1}^T\rE\Big\{\Big[\bbF_j\bgma_{j+\tau}^*\bbB_2\Big(\ell\bbI-\bbB_2^*\bbM\bbB_2\Big)^{-1}\bbB_2^*\bgma_{j+\tau}\bbF_j^*+\\
&&\bbF_{j+\tau}\bgma_j^*\bbB_2\Big(\ell\bbI-\bbB_2^*\bbM\bbB_2\Big)^{-1}\bbB_2^*\bgma_j\bbF_{j+\tau}^*\Big]\Big|\bgma_1,\cdots,\bgma_{T+\tau}\Big\}\\
&=&\frac1{2T}\rE\rtr\Big[\Big(\ell\bbI-\bbB_2^*\bbM\bbB_2\Big)^{-1}\bbB_2^*\sum_{j=1}^T(\bgma_{j+\tau}\bgma_{j+\tau}^*+\bgma_j\bgma_j^*)\bbB_2\Big]\bbI_{k(q+1)}\\
&=&\frac1T\sum_{j=1}^T\rE\bgma_j^*\bbB_2\Big(\ell\bbI-\bbB_2^*\bbM\bbB_2\Big)^{-1}\bbB_2^*\bgma_j\bbI_{k(q+1)}\\
&=&C_0\bbI_{k(q+1)}+o(\textbf{1}).
\eqn
Similarly, we have
\bqn
\rE(P_2)&=&\frac12C_1\bbH(\tau)+o(\textbf{1}),\\
\rE(P_3)&=&C_1\bbH_L(\tau)+\frac12(C_0\bbH_L(2\tau)+C_2\bbI_{k(q+1)})+o(\textbf{1}),\\
\rE(P_4)&=&C_1\bbH_U(\tau)+\frac12(C_0\bbH_U(2\tau)+C_2\bbI_{k(q+1)})+o(\textbf{1}).\\
\eqn
Here $\bbH_L(\tau)$ and $\bbH_U(\tau)$ denote the lower and upper part of $\bbH(\tau)$ with the rest entries being 0. Furthermore, denote $\bbH_L(0)\equiv\bbH_U(0)\equiv\bbI_{k(q+1)}$ and hence $\bbH(0)=\bbH_L(0)+\bbH_U(0)=2\bbI_{k(q+1)}$. Consider $i=j\pm p\tau$ for $p=0,1,\cdots,\big[\frac q\tau\big]$.\\
When $p=0$, we have $\rE(\bbP1)+\rE(\bbP2)=\frac12C_0\bbH(0)+\frac12C_1\bbH(\tau)+o(\textbf{1}).$\\
When $p=1$, we have $\rE(\bbP3)+\rE(\bbP4)=\frac12C_2\bbH(0)+C_1\bbH(\tau)+\frac12C_0\bbH(2\tau)+o(\textbf{1}).$\\
When $p=2$, we have part of $\rE(\bbP5)$ is $\frac12C_3\bbH(\tau)+C_2\bbH(2\tau)+\frac12C_1\bbH(3\tau)+o(\textbf{1}).$\\
$\vdots$\\
When $p=\big[\frac q\tau\big]$, we have part of $\rE(\bbP5)$ is
\bqn
&&\frac12C_{\big[\frac q\tau\big]+1}\bbH(\big[\frac q\tau\big]\tau-\tau)+C_{\big[\frac q\tau\big]}\bbH(\big[\frac q\tau\big]\tau)+\frac12C_{\big[\frac q\tau\big]-1}\bbH(\big[\frac q\tau\big]\tau+\tau)+o(\textbf{1})\\
&=&\frac12C_{\big[\frac q\tau\big]+1}\bbH(\big[\frac q\tau\big]\tau-\tau)+C_{\big[\frac q\tau\big]}\bbH(\big[\frac q\tau\big]\tau)+o(\textbf{1}).
\eqn
When $p=\big[\frac q\tau\big]+1$, we have part of $\rE(\bbP5)$ is
\bqn
&&\frac12C_{\big[\frac q\tau\big]}\bbH(\big[\frac q\tau\big]\tau+2\tau)+C_{\big[\frac q\tau\big]+1}\bbH(\big[\frac q\tau\big]\tau+\tau)+\frac12C_{\big[\frac q\tau\big]+2}\bbH(\big[\frac q\tau\big]\tau)+o(\textbf{1})\\
&=&\frac12C_{\big[\frac q\tau\big]+2}\bbH(\big[\frac q\tau\big]\tau)+o(\textbf{1}).
\eqn
Next, we want to show that $\bbP_i\to\rE(\bbP_i)$ a.s. Since all the $\bbP_i$'s are of finite dimension, it suffices to show the a.s convergence entry-wise. Denote the $(u,v)$-entry of a matrix $\bbA$ by $\bbA_{(u,v)}$. For $i=1$, define
$\alpha_j=\bgma_j^*\bbB_2\Big(\ell\bbI-\bbB_2^*\bbM\bbB_2\Big)^{-1}\bbB_2^*\bgma_j$. Then for any positive integer $s$, applying Lemma \ref{Burkholder}, we have
\bqa\label{P1}
&&\rE(|\bbP_{1(i1,i2)}-\frac{\alpha_j}2\delta_{(i_1,i_2)}-\frac{\alpha_{j+\tau}}2\delta_{(i_1,i_2)}|^{2s})\non
&=&\rE\Big\{\frac1{2T}\sum_{j=1}^T\Big[\alpha_{j+\tau}(\bbF_j\bbF_j^*)_{(i_1,i_2)}+\alpha_j(\bbF_{j+\tau}\bbF_{j+\tau}^*)_{(i_1,i_2)}\Big]-\frac{\alpha_j}2\delta_{(i_1,i_2)}-\frac{\alpha_{j+\tau}}2\delta_{(i_1,i_2)}\Big\}^{2s}\non
&\le&2^{2s-1}\rE\Big[\frac1{2T}\sum_{j=1}^T\alpha_{j+\tau}(\bbF_j\bbF_j^*)_{(i_1,i_2)}-\frac{\alpha_{j+\tau}}2\delta_{(i_1,i_2)}\Big]^{2s}+\non
&&2^{2s-1}\rE\Big[\frac1{2T}\sum_{j=1}^T\alpha_j(\bbF_{j+\tau}\bbF_{j+\tau}^*)_{(i_1,i_2)}-\frac{\alpha_j}2\delta_{(i_1,i_2)}\Big]^{2s}\non
&=&\frac12\rE\Big[\frac1T\sum_{j=1}^T\alpha_{j+\tau}(\bbF_j\bbF_j^*)_{(i_1,i_2)}-\alpha_{j+\tau}\delta_{(i_1,i_2)}\Big]^{2s}+\non
&&\frac12\rE\Big[\frac1T\sum_{j=1}^T\alpha_j(\bbF_{j+\tau}\bbF_{j+\tau}^*)_{(i_1,i_2)}-\alpha_j\delta_{(i_1,i_2)}\Big]^{2s}\non
&=&\frac12\rE\Big\{\frac1T\sum_{j=1}^T\bbF_j^*\Big[\alpha_{j+\tau}^{(i_1,i_2)}\Big]\bbF_j-\rtr\Big[\alpha_{j+\tau}^{(i_1,i_2)}\Big]\Big\}^{2s}+\non
&&\frac12\rE\Big\{\frac1T\sum_{j=1}^T\bbF_{j+\tau}^*\Big[\alpha_j^{(i_1,i_2)}\Big]\bbF_{j+\tau}-\rtr\Big[\alpha_j^{(i_1,i_2)}\Big]\Big\}^{2s}.
\eqa
Here $\big[a^{(u,v)}\big]$ denotes the matrix with the $(u,v)$-entry being $a$ and 0 elsewhere.
By the truncation of $\ep_{ij}$ and the fact that $\|(\ell\bbI-\bbB_2^*\bbM\bbB_2)^{-1}\|_o\le \eta^{-1}$ with $\eta=\ell-d_c>0$, both $|\alpha_j|$ and $|\alpha_{j+\tau}|$ are bounded from above, say, by $C$. Also notice that $|\bbF^{\tau}_{ij}|< T^{1/4}$ and $\rE |\bbF^{\tau}_{ij}|^4<M$. Similar to the proof of Lemma 9.1 in Bai and Silverstein (2010), we have
\bqn
&&\rE\Big\{\frac1T\sum_{j=1}^T\bbF_j^*\Big[\alpha_{j+\tau}^{(i_1,i_2)}\Big]\bbF_j-\rtr\Big[\alpha_{j+\tau}^{(i_1,i_2)}\Big]\Big\}^{2s}
\le \frac{C^{2s}}{T^s}\sum_{l=1}^s \big(M^ll^{2s}+l^{4s}\big)\\
&&\rE\Big\{\frac1T\sum_{j=1}^T\bbF_{j+\tau}^*\Big[\alpha_j^{(i_1,i_2)}\Big]\bbF_{j+\tau}-\rtr\Big[\alpha_j^{(i_1,i_2)}\Big]\Big\}^{2s}
\le \frac{C^{2s}}{T^s}\sum_{l=1}^s \big(M^ll^{2s}+l^{4s}\big).
\eqn
Substituting the above back to (\ref{P1}) and choosing $s\ge 2$, we have $$\bbP_{1(i_1,i_2)}-\frac{\alpha_j}2\delta_{(i_1,i_2)}-\frac{\alpha_{j+\tau}}2\delta_{(i_1,i_2)}=o_{a.s}(1).$$ Again, by the almost sure and uniform
convergence of $\alpha_j$ and $\alpha_{j+\tau}$ to $C_0$, we have $$\frac{\alpha_j}2\delta_{(i_1,i_2)}+\frac{\alpha_{j+\tau}}2\delta_{(i_1,i_2)}-C_0\delta_{(i_1,i_2)}=o_{a.s}(1).$$ Therefore, we have shown that $P_1-\rE(P_1)=o_{a.s.}(1)$. Results for $i=2,3,4,5$ can be shown in a similar way.

Denote $\alpha=\frac{-\frac{cm}2}{1-\frac{c^2m^2}{2x_1}}$ and $\beta=\frac{-\frac{cm}2}{1-\frac{c^2m^2}{4x_1}}$, then we have $C_p=\alpha\beta^p$. Note that $\bbH(p\tau)=\textbf{0}$ for $p>[q/\tau]$, and we have, with probability 1 that
\bqn
&&\bbW\bbB_2\Big(\ell\bbI-\bbB_2^*\bbM\bbB_2\Big)^{-1}\bbB_2^*\bbW^*\\
&\to&(\frac12C_0+\frac12C_2)\bbH(0)+(\frac32C_1+\frac12C_3)\bbH(\tau)+\sum_{p=2}^{\infty}(\frac12C_{p-2}+C_p+\frac12C_{p+2})\bbH(p\tau)\\
&=&\frac{\alpha}2\Big[(1+\beta^2)\bbH(0)+(3\beta+\beta^3)\bbH(\tau)+(1+\beta^2)^2\sum_{p=2}^{\infty}\beta^{p-2}\bbH(p\tau)\Big]\\
&=&\frac{\alpha}2\Big[(1+\beta^2)\bbH(0)+(3\beta+\beta^3)\bbH(\tau)+(1+\beta^2)^2\sum_{p=2}^{\infty}\beta^{p-2}\big(\bbH_L(p\tau)+\bbH_U(p\tau)\big)\Big]\\
&=&\frac{\alpha}2\Big[(1+\beta^2)\bbH(0)+(3\beta+\beta^3)\bbH(\tau)+(1+\beta^2)^2\sum_{p=2}^{\infty}\beta^{p-2}\big(\bbJ_L(p\tau)+\bbJ_U(p\tau)\big)\otimes\bbI_k\Big]\\
&=&\frac{\alpha}2\Big[(1+\beta^2)\bbH(0)+(3\beta+\beta^3)\bbH(\tau)+(1+\beta^2)^2\sum_{p=2}^{\infty}\beta^{p-2}\big(\bbJ_L^p(\tau)+\bbJ_U^p(\tau)\big)\otimes\bbI_k\Big]\\
&=&\frac{\alpha}2\Big[(1+\beta^2)\bbH(0)+(3\beta+\beta^3)\bbH(\tau)+(1+\beta^2)^2\Big(\bbJ_L^2(\tau)(\bbI-\beta\bbJ_L(\tau))^{-1}\\
&&+\bbJ_U^2(\tau)(\bbI-\beta\bbJ_U(\tau))^{-1}\Big)\otimes\bbI_k\Big].
\eqn
Note that
\bqn
\bbJ_L(\tau)(\bbI-\beta\bbJ_L(\tau))^{-1}&=&\frac1\beta[\bbI-(\bbI-\beta\bbJ_L(\tau))](\bbI-\beta\bbJ_L(\tau))^{-1}\\
&=&\frac1\beta[(\bbI-\beta\bbJ_L(\tau))^{-1}-\bbI]\\
\bbJ_L^2(\tau)(\bbI-\beta\bbJ_L(\tau))^{-1}&=&\frac1\beta\bbJ_L(\tau)[(\bbI-\beta\bbJ_L(\tau))^{-1}-\bbI]\\
&=&\frac1{\beta^2}[(\bbI-\beta\bbJ_L(\tau))^{-1}-\bbI]-\frac1\beta\bbJ_L(\tau).
\eqn
Similarly,
\bqn
\bbJ_U^2(\tau)(\bbI-\beta\bbJ_U(\tau))^{-1}&=&\frac1\beta\bbJ_U(\tau)[(\bbI-\beta\bbJ_U(\tau))^{-1}-\bbI]\\
&=&\frac1{\beta^2}[(\bbI-\beta\bbJ_U(\tau))^{-1}-\bbI]-\frac1\beta\bbJ_U(\tau).
\eqn
Therefore, we have
\bqn
&&(1+\beta^2)^2\Big(\bbJ_L^2(\tau)(\bbI-\beta\bbJ_L(\tau))^{-1}+\bbJ_U^2(\tau)(\bbI-\beta\bbJ_U(\tau))^{-1}\Big)\otimes\bbI_k\\
&=&(1+\beta^2)^2\Big(\frac1{\beta^2}(\bbI-\beta\bbJ_L(\tau))^{-1}+\frac1{\beta^2}(\bbI-\beta\bbJ_U(\tau))^{-1}\\
&&-\frac2{\beta^2}\bbI-\frac1\beta\big(\bbJ_L(\tau)+\bbJ_U(\tau)\big)\Big)\otimes\bbI_k\\
&=&\frac{(1+\beta^2)^2}{\beta^2}\Big[(\bbI-\beta\bbJ_L(\tau))^{-1}\Big(2\bbI-\beta\bbJ_L(\tau)-\beta\bbJ_U(\tau)\Big)(\bbI-\beta\bbJ_U(\tau))^{-1}\Big]\otimes\bbI_k\\
&&-\frac{(1+\beta^2)^2}{\beta^2}\bbH(0)-\frac{(1+\beta^2)^2}{\beta}\bbH(\tau)\\
&=:&\frac{(1+\beta^2)^2}{\beta^2}\bbG(\tau)-\frac{(1+\beta^2)^2}{\beta^2}\bbH(0)-\frac{(1+\beta^2)^2}{\beta}\bbH(\tau).
\eqn
Hence, we have a.s.
\bqn
&&\bbB_1^*\bgL\bbW\bbB_2\Big(\ell\bbI-\bbB_2^*\bbM\bbB_2\Big)^{-1}\bbB_2^*\bbW^*\bgL^*\bbB_1\\
&\to&\frac\alpha2\bbQ^{1/2}\Big[\Big(1+\beta^2-\frac{(1+\beta^2)^2}{\beta^2}\Big)\bbH(0)+\Big(3\beta+\beta^3-\frac{(1+\beta^2)^2}{\beta}\Big)\bbH(\tau)\\
&&+\frac{(1+\beta^2)^2}{\beta^2}\bbG(\tau)\Big]\bbQ^{1/2}.
\eqn
Last, we want to show that with probability 1,
\bqn
\bbB_1^*\bgL\bbW\bbB_2\Big(\ell\bbI-\bbB_2^*\bbM\bbB_2\Big)^{-1}\bbB_2^*\bbM\bbB_1\to\textbf{0}
\eqn
and
\bqn
\bbB_1^*\bbM\bbB_2\Big(\ell\bbI-\bbB_2^*\bbM\bbB_2\Big)^{-1}\bbB_2^*\bbW^*\bgL^*\bbB_1\to\textbf{0}.
\eqn
Note that
\bqn
\bbB_1^*\bgL \bbW\bbB_2\Big(\ell\bbI-\bbB_2^*\bbM\bbB_2\Big)^{-1}\bbB_2^*\bbM\bbB_1
=\bbB_1^*\bgL \bbW\Big(\bbI-\frac1\ell\bbB_2\bbB_2^*\bbM\Big)^{-1}\bbB_1-\bbB_1^*\bgL\bbW\bbB_1
\eqn
and that
\bqn
\bbB_1^*\bbM\bbB_2\Big(\ell\bbI-\bbB_2^*\bbM\bbB_2\Big)^{-1}\bbB_2^*\bgL^*\bbW\bbB_1
=\bbB_1^*\Big(\bbI-\frac1\ell\bbM\bbB_2\bbB_2^*\Big)^{-1}\bbW^*\bgL^*\bbB_1-\bbB_1^*\bbW^*\bgL^*\bbB_1.
\eqn
Hence, by $\bbW\bbB_1=o_{a.s.}(\textbf{1})$ and $\bbB_1^*\bbW^*=o_{a.s.}(\textbf{1})$, it suffices to show with probability 1 that,
\bqn
\bbB_1^*\bgL \bbW\Big(\bbI-\frac1\ell\bbB_2\bbB_2^*\bbM\Big)^{-1}\bbB_1\to\textbf{0}
\eqn
and
\bqn
\bbB_1^*\Big(\bbI-\frac1\ell\bbM\bbB_2\bbB_2^*\Big)^{-1}\bbW^*\bgL^*\bbB_1\to\textbf{0}.
\eqn
By $\bbB_1=\bgL(\bgL^*\bgL)^{-1/2}$, we have
\bqn
&&\bbB_1^*\bgL \bbW\Big(\bbI-\frac1\ell\bbB_2\bbB_2^*\bbM\Big)^{-1}\bbB_1\\
&=&(\bgL^*\bgL)^{1/2}\bbW\Big(\bbI-\frac1\ell\bbB_2\bbB_2^*\bbM\Big)^{-1}\bgL(\bgL^*\bgL)^{-1/2}\\
&=&\bbQ^{1/2}\bbW\Big(\bbI-\frac1\ell\bbB_2\bbB_2^*\bbM\Big)^{-1}\bgL\bbQ^{-1/2}.
\eqn
By law of large numbers, we have with probability 1 that,
\bqn
&&\bbW\Big(\bbI-\frac1\ell\bbB_2\bbB_2^*\bbM\Big)^{-1}\bgL\\
&=&\frac1{2T}\sum_{j=1}^T(\bbF_j\bbep_{j+\tau}^*+\bbF_{j+\tau}\ep_j^*)\Big(\bbI-\frac1\ell\bbB_2\bbB_2^*\bbM\Big)^{-1}\bgL\\
&\to&\rE(\bbF_1\bbep_{1+\tau}^*+\bbF_{1+\tau}\bbep_1^*)\Big(\bbI-\frac1\ell\bbB_2\bbB_2^*\bbM\Big)^{-1}\bgL\\
&=&\rE\rE\Big((\bbF_1\bbep_{1+\tau}^*+\bbF_{1+\tau}\bbep_1^*)\Big(\bbI-\frac1\ell\bbB_2\bbB_2^*\bbM\Big)^{-1}\bgL\Big|\bbep_1,\cdots,\bbep_{T+\tau}\Big)\\
&=&o(\textbf{1}).
\eqn
Hence, we have with probability 1
\bqn
\bbB_1^*\bgL \bbW\Big(\bbI-\frac1\ell\bbB_2\bbB_2^*\bbM\Big)^{-1}\bbB_1
=\bbQ^{1/2}\bbW\Big(\bbI-\frac1\ell\bbB_2\bbB_2^*\bbM\Big)^{-1}\bgL\bbQ^{-1/2}
=o(\textbf{1}).
\eqn
Similarly,
\bqn
\bbB_1^*\Big(\bbI-\frac1\ell\bbM\bbB_2\bbB_2^*\Big)^{-1}\bbW^*\bgL^*\bbB_1=o_{a.s.}(\textbf{1}).
\eqn
Therefore, $\ell$ should satisfy
\bqn
&&\det\Big|\bbQ^{1/2}\bbH(\tau)\bbQ^{1/2}+\Big(\bbB_1^*\big(\bbM-\ell\bbI\big)^{-1}\bbB_1\Big)^{-1}+\\
&&\frac\alpha2\bbQ^{1/2}\Big[\Big(1+\beta^2-\frac{(1+\beta^2)^2}{\beta^2}\Big)\bbH(0)+\Big(3\beta+\beta^3-\frac{(1+\beta^2)^2}{\beta}\Big)\bbH(\tau)\\
&&+\frac{(1+\beta^2)^2}{\beta^2}\bbG(\tau)\Big]\bbQ^{1/2}\Big|\to 0.
\eqn
Recall $\bbB_1=\bgL(\bgL^*\bgL)^{-1/2}$. Our next goal is to find the limit of
\bqn
\bbB_1^*\big(\bbM-\ell\bbI\big)^{-1}\bbB_1.
\eqn
Define $\bbA=\bbM-\ell\bbI$ and $\bbA_k=\bbA-(\bgma_{k+\tau}+\bgma_{k-\tau})\bgma_k^*-\bgma_k(\bgma_{k+\tau}+\bgma_{k-\tau})^*$, then we have the following lemmas, with proofs given in the Appendix.
\begin{lemma}\label{gAx}Let $\bbx\in \mathbb{C}_1^n:=\{\bbx\in\mathbb{C}^n: \|\bbx\|=1\}$ be given. For $r\ge1$, we have
\bqn
\rE|\bgma_k^*\bbA_k^{-1}\bbx|^{2r}\le KT^{-r}
\eqn
for some $K>0$.
\end{lemma}
\begin{lemma}\label{xAy} For any $\bbx,\bby\in \mathbb{C}_1^n$, we have $\bbx^*\bbA^{-1}\bby\to-\frac{\bbx^*\bby}{\frac{cm}{1-c^2m^2+\sqrt{1-c^2m^2}}+\ell}$ a.s.
\end{lemma}
Finally, we have
\bqn
&&\det\Big|\bbQ^{1/2}\bbH(\tau)\bbQ^{1/2}-\Big(\frac{cm(\ell)}{1-c^2m^2(\ell)+\sqrt{1-c^2m^2(\ell)}}+\ell\Big)\bbI_{k(q+1)}\\
&&+\frac\alpha2\bbQ^{1/2}\Big[\Big(1+\beta^2-\frac{(1+\beta^2)^2}{\beta^2}\Big)\bbH(0)+\Big(3\beta+\beta^3-\frac{(1+\beta^2)^2}{\beta}\Big)\bbH(\tau)\\
&&+\frac{(1+\beta^2)^2}{\beta^2}\bbG(\tau)\Big]\bbQ^{1/2}\Big|=0,
\eqn
or equivalently
\bqa
&&\det\Big|\frac\alpha2\Big[\Big(1+\beta^2-\frac{(1+\beta^2)^2}{\beta^2}\Big)\bbH(0)+\Big(3\beta+\beta^3-\frac{(1+\beta^2)^2}{\beta}+\frac2\alpha\Big)\bbH(\tau)\non
&&+\frac{(1+\beta^2)^2}{\beta^2}\bbG(\tau)\Big]-\Big(\frac{cm(\ell)}{1-c^2m^2(\ell)+\sqrt{1-c^2m^2(\ell)}}+\ell\Big)\bbQ^{-1}\Big|=0.
\label{det}
\eqa
When $\tau>q$, one has $\bbH(\tau)=\textbf{0}$, $\bbG(\tau)=2\bbI$ and (\ref{det}) reduces to
\bqn
\det\Big|\alpha(1+\beta^2)\bbI-\Big(\frac{cm(\ell)}{1-c^2m^2(\ell)+\sqrt{1-c^2m^2(\ell)}}+\ell\Big)\bbQ^{-1}\Big|=0,
\eqn
or
\bqn
\det\Big|\frac{cm(\ell)}{1-c^2m^2(\ell)+\sqrt{1-c^2m^2(\ell)}}(\bbQ+\bbI)+\ell\bbI\Big|=0.
\eqn
Let $\lambda$ be an eigenvalue of $\bbQ$, then we have
\bqa\label{det1}
\frac{cm(\ell)}{1-c^2m^2(\ell)+\sqrt{1-c^2m^2(\ell)}}\big(1+\frac1{\la}\big)+\frac\ell{\la}=0.
\eqa
When $\tau=q$, (\ref{det}) reduces to
\bqn
\det\Big|\alpha(1+\beta^2)\bbI+\Big(1+\frac\alpha2(3\beta+\beta^3)\Big)\bbH(q)-\Big(\frac{cm(\ell)}{1-c^2m^2(\ell)+\sqrt{1-c^2m^2(\ell)}}+\ell\Big)\bbQ^{-1}\Big|=0.
\eqn
Writing
$$\bbH(q)=\begin{pmatrix}
  0 & \cdots & \cdots & 0 & 1 \\
  \vdots & \ddots & \ddots &  \ddots & 0\\
 \vdots &  \ddots &\ddots&\ddots & \vdots  \\
  0 & \ddots& \ddots & \ddots & \vdots   \\
  1 &0 & \cdots& \cdots & 0  \\
 \end{pmatrix}\otimes\bbI_k,$$
 one can easily verify that the eigenvalues of $\bbH(q)$ are 1, $-1$ and 0, with multiplicity $k$, $k$ and $k(q-1)$, respectively.\\
Suppose that $\bbH(q)$ and $\bbQ$ are commutative,  that is, there is a common orthogonal matrix
$\bbO$ simultaneously diagnalizing the two matrices, i.e., we have $\bbH(q)=\bbO\bbD^{\bbH}\bbO'$ and $\bbQ=\bbO\bbD^{\bbQ}\bbO'$,
where $\bbD^{\bbH}= \mbox{diag}[a_1,\cdots,a_{k(q+1)}]$ and $\bbD^{\bbQ}= \mbox{diag}[\la_1,\cdots,\la_{k(q+1)}]$.
Then, (\ref{det}) further reduces to
\bqn
&&\Big(1+\frac\alpha2(3\beta+\beta^3)\Big)a_j\\
&=&\frac{cm(\ell)}{1-c^2m^2(\ell)+\sqrt{1-c^2m^2(\ell)}}\big(1+\frac1{\la_j}\big)+\frac\ell{\la_j},\quad j=1,\cdots,k(q+1).
\eqn
Substituting $\alpha=\frac{-\frac{cm}2}{1-\frac{c^2m^2}{2x_1}}$ and $\beta=\frac{-\frac{cm}2}{1-\frac{c^2m^2}{4x_1}}$, for $j=1,\cdots,k(q+1)$, we have
\bqa\label{eqaj}
a_j&=&\Big(\frac12+\frac1{1-c^2m^2(\ell)+\sqrt{1-c^2m^2(\ell)}}\Big)^{-1}\times\non
&&\Big[\frac{cm(\ell)}{1-c^2m^2(\ell)+\sqrt{1-c^2m^2(\ell)}}\big(1+\frac1{\la_j}\big)+\frac\ell{\la_j}\Big]\non
&=:&g_j(\ell).
\eqa
Notice that (\ref{det1}) is a special case of (\ref{eqaj}) for $a_j=0$.\\
Note that when $x$ is outside the support $[-d(c),d(c)]$ of the LSD of $\bbM_n(\tau)$, $m_2(x)\ne0$. Hence, we have
\bqn
&cm_1(x)((1-c-cxm_1(x))^2-c^2x^2m_2^2(x))&\\
&+x(1-c^2m_1^2(x)+c^2m_2^2(x))(1-c-cxm_1(x))=0.&
\eqn
Let $x\downarrow d(c):=d$ and we have $m_2(x)\to 0$ and $m_1(x)\to m_1(d)$ satisfying
\bqn
&&cm_1(d)(1-c-cdm_1(d))^2+d(1-c^2m_1^2(d))(1-c-cdm_1(d))\\
&=&[1-c-cdm_1(d)][cm_1(d)(1-c-cdm_1(d))+d(1-c^2m_1^2(d))]=0,
\eqn
from which we have $m_1(d)=\frac{1-c-\sqrt{(1-c)^2+8d^2}}{4cd}$.\\
Rewrite
\bqn
g_j(\ell)&=&\frac{cm(\ell)}{\frac32-\frac12c^2m^2(\ell)+\frac12\sqrt{1-c^2m^2(\ell)}}\big(1+\frac1{\la_j}\big)\\
&&+\Big(\frac12+\frac1{1-c^2m^2(\ell)+\sqrt{1-c^2m^2(\ell)}}\Big)^{-1}\frac\ell{\la_j}\\
&:=&g_{j1}(\ell)+g_{j2}(\ell).
\eqn
We will show that $g_j(\ell)$ is increasing over $(d(c),\infty)$ by showing so are $g_{j1}$ and $g_{j2}$.
By definition, over $(d(c),\infty)$, $m$ is an increasing function taking negative values and $m^2$ is a decreasing function taking positive values. Hence, it is easy to see that $g_{j1}(\ell)$ is increasing over $(d(c),\infty)$. For $g_{j2}(\ell)$, define $h(\ell)=1-c^2m^2(\ell)+\sqrt{1-c^2m^2(\ell)}$ and rewrite $g_{j2}(\ell)=\frac{\ell h(\ell)}{\la_j(1+\frac{h(\ell)}2)}$.
It is easy to see that $h(\ell)>0$ and $h'(\ell)>0$ over $(d(c),\infty)$. Hence we have
\bqn
g_{j2}'(\ell)&=&\frac{[h(\ell)+\ell h'(\ell)](1+\frac{h(\ell)}2)-\ell h(\ell)\frac{h'(\ell)}2}{\la_j(1+\frac{h(\ell)}2)^2}
=\frac{h(\ell)+\ell h'(\ell)+\frac{h^2(\ell)}2}{\la_j(1+\frac{h(\ell)}2)^2}>0.
\eqn
By symmetry, $g_j(\ell)$ is increasing over $(-\infty,-d(c))$ as well. Therefore, based on the sign of $g_j(d(c))$, we have the following cases to consider.

Case I. $g_j(d(c))\ge0$:\\
i. If $a_j>g_j(d(c))$, then (\ref{eqaj}) has one solution in $(d(c),\infty)$ and no solution in $(-\infty,-d(c))$.\\
ii. If $g(d(c))\ge a_j\ge-g_j(d(c))=g_j(-d(c))$, then (\ref{eqaj}) has no solution in $(d(c),\infty)$ and $(-\infty,-d(c))$.\\
iii. If $a_j<-g_j(d(c))$, then (\ref{eqaj}) has one solution in $(-\infty,-d(c))$ and no solution in $(d(c),\infty)$.

Case II. $g_j(d(c))<0$:\\
i. If $a_j\ge -g_j(d(c))$, then (\ref{eqaj}) has one solution in $(d(c),\infty)$ and no solution in $(-\infty,-d(c))$.\\
ii. If $g(d(c))< a_j <-g_j(d(c))=g_j(-d(c))$, then (\ref{eqaj}) has one solution in $(d(c),\infty)$ and one solution in $(-\infty,-d(c))$.\\
iii. If $a_j\le g_j(d(c))$, then (\ref{eqaj}) has one solution in $(-\infty,-d(c))$ and no solution in $(d(c),\infty)$.

\begin{remark} In real application, compared with the noise component, the loading matrix $\bgL$ dominates. As a result, all the eigenvalues of $\bbQ=\bgL^*\bgL$ are large (more precisely, they are of the same order as $n$). Hence, we can assume that $\bbQ^{-1}=\textbf{0}$. Thus the commutative assumption of $\bbH(q)$ and $\bbQ$ can be relaxed. Moreover, under this case, we always have $g_j(d(c))<0$.
\end{remark}
\begin{remark}
For the same reason as stated before Remark 4.1, $d_c$ is replaced by $(1+an^b)d_c$ in practice. Simulation indicates a fit of $a=0.1, b=-1/3$.
\end{remark}

Notice that all the eigenvalues of $\bbH(q+1)$ are 0, while for $\bbH(q)$, $k$ eigenvalues are 1 and $k$ eigenvalues are $-1$, with the rest being 0. Making use of such difference and applying the above analysis to the cases that $\tau=q$ and $\tau=q+1$ gives an estimate of $q$. Together with the estimation of $k(q+1)$, we easily obtain the estimate of $k$. A numerical demonstration is given in the simulation.

\section{Estimate of $\sigma^2$}\label{sigma}
\setcounter{equation}{0}
\def\theequation{\thesection.\arabic{equation}}
\setcounter{subsection}{0}
The above estimation is based on the assumption that $\sigma^2$, the variance of the noise part is given. In practice, it is often the case that $\sigma^2$ is unknown. To this end, we can estimate $\sigma^2$ by employing the properties of the MP law. More precisely, we first estimate the left boundary of the support of the MP law by the smallest sample eigenvalue of $\bbPhi_n(0)$, say $\hat\la_{1}$ (the eigenvalues are arranged in ascending order), and estimate the right boundary by $\frac{(1+\sqrt c)^2}{(1-\sqrt c)^2}\hat\la_{1}$. An iteration is then applied.  The initial estimator of $\sigma^2$, say $\hat\sigma^2_{(0)}$, is obtained as the sample mean of the sample eigenvalues of $\bbPhi_n(0)$ that lie within the interval $[\hat\la_{m+1}, \hat\la_{n-m}]$, where $m\ge 0$ is such that $\hat\la_{n-m}$ is the largest eigenvalue of $\bbPhi_n(0)$ less than $ \frac{(1+\sqrt c)^2(1+2n^{-2/3})}{(1-\sqrt c)^2}\hat\la_{1}$. For $i\ge1$, we obtain the updated estimator $\hat\sigma^2_{(i)}$ by taking the sample mean of the sample eigenvalues of $\bbPhi_n(0)$ that lie within the interval $[(1-\sqrt c)^2\hat\sigma^2_{(i-1)}, (1+\sqrt c)^2(1+2n^{-2/3})\hat\sigma^2_{(i-1)}]$. The iteration stops once we have $\hat\sigma^2_{(\ell-1)}=\hat\sigma^2_{(\ell)}$ for some $\ell$ and our estimator $\hat\sigma^2:=\hat\sigma^2_{(\ell)}$.  As shown in the simulation, our estimation of $k$ and $q$ still works well with such estimator.
\section{Simulation}
\setcounter{equation}{0}
\def\theequation{\thesection.\arabic{equation}}
\setcounter{subsection}{0}
Table 1 presents a simulation about the result discussed above, displaying the largest 13 absolute values of the eigenvalues
for lags $\tau$ from 0 to 5. Here
\begin{equation}
\bbR_{t}=\sum_{i=0}^{q}\bgL_{i}\bbf_{t-i}+\bbe_t,\quad t=1,...,T\label{eq:}\end{equation}
 where $\bbf_{t}$ 's are factors of length $k$; $\bgL_{i}, i=0,...,q$
is a constant time-invariant matrix of size $n\times k$, $\bbe_{t}$
is the error term and $q$ is the lag of the model. In addition, assume that:
$\bbe_{t}$ are i.i.d. random variables with $\bbe_t\sim N\left(\mathbf{0},\sigma^{2}\bbI_{\mathit{n}}\right)$ and $\bbf_{t}$ are i.i.d. random variables with $\bbf_{t}\sim N\left(\mathbf{0},\sigma_{\mathit{f}}^{2}\bbI_{\mathit{k}}\right)$, independent of $\bbe_{t}$. $\bgL_{i}=\left[\begin{array}{cccccc}
\bgL_{1}^{i} & \bgL_{2}^{i} & . & . & . & \bgL_{k}^{i}\end{array}\right]$ where $\bgL_{j}^{i}\quad i=0,...,q;j=1,...,k$ is a vector of
length $n$ and is given by $\bgL_{j}^{i}=\beta\mathbf{1}_{n}+\ep_{ij}$ where $\mathbf{1}_{n}$ is a vector of 1's and $\ep_{ij}$ are i.i.d.
random variables with $\ep_{ij}\sim N\left(\mathbf{0},\sigma_{\ep}^{2}\bbI_{\mathit{n}}\right)$.
For $n=450$, $T=500$, $k=2$, $q=2$, $\beta=1.0$, $\sigma_{f}^{2}=4$, $\sigma^2=1$ and $\sigma_{\ep}^{2}=0.25$, we have $c=0.9$, $b_c=(1+\sqrt c)^2=3.7974$ and $d_c=1.8573$. 
Eigenvalues of $\bbQ$ are from 95 to 285, making $\bbQ^{-1}\sim \textbf{0}$.

When $\tau=0$, and $\sigma^2 = 1$ is known, using the phase transition point $b_c=(1+\sqrt c)^2 \sigma^2=3.7974$, we see that the number of spotted spikes is 6, which estimates $k(q+1)$. When for $\tau=q+1$, we have $\bbH(\tau)=\textbf{0}$. Moreover, as $\bbQ^{-1}\sim \textbf{0}$, we have $g_j(d(c))\sim -0.4284<0$. That is, our Case II (ii) applies for all the $k(q+1)$ 0 eigenvalues of $\bbH(q+1)$, making the number of spikes $2k(q+1)$ as verified by applying the phase transition point $d_c=1.8573$. For $\tau=q$, $\bbH(\tau)$ has $k$ eigenvalues of 1, $k$ eigenvalues of $-1$ and $k(q-1)$ eigenvalues of $0$ with Case II (i),(iii) and (ii) applicable, respectively. Thus, we have $k+k+2k(q-1)=2kq<2k(q+1)$ eigenvalues in this case. Again, this agrees with the use of the phase transition point $d_c=1.8573$. In other words, the number of spikes first jumps to $2k(q+1)$ at $\tau=3$ which estimates $q+1$. The estimation of $k$ is obvious.

When $\sigma^2=1$ is unknown, using technique as in Section \ref{sigma}, one has $\hat\sigma^2=0.9894$. It then follows that $\hat b_c=(1+\sqrt c)^2(1+2n^{-2/3})\hat\sigma^2=3.8851$, $\hat d_c=1.8616$ (with rescale factor $1+0.1n^{-1/3}$), which gives the same estimates as above.
\begin{table}[]
\centering{}\begin{tabular}{|c|c|c|c|c|c|}
\hline
$\tau=0$ & \multicolumn{1}{c|}{$\tau=1$} & $\tau=2$ & $\tau=3$ & $\tau=4$ & $\tau=5$\tabularnewline
\hline
\hline
\textcolor[rgb]{0.501961,0,0}{\textbf{10031.2366}}  & \textcolor[rgb]{0.501961,0,0}{\textbf{6227.5906}}  & \textcolor[rgb]{0.501961,0,0}{\textbf{2865.1554}} & \textcolor[rgb]{0.501961,0,0}{\textbf{640.5761}} & \textbf{155.9377} & \textbf{128.6870}\tabularnewline
\hline
\textcolor[rgb]{0.501961,0,0}{\textbf{534.5839}}  & \textcolor[rgb]{0.501961,0,0}{\textbf{363.7782}} & \textcolor[rgb]{0.501961,0,0}{\textbf{258.4859}} & \textcolor[rgb]{0.501961,0,0}{\textbf{48.9667}} & \textbf{46.5224}  & \textbf{92.3756}\tabularnewline
\hline
\textcolor[rgb]{0.501961,0,0}{\textbf{473.1639}}  & \textcolor[rgb]{0.501961,0,0}{\textbf{325.8391}}  & \textcolor[rgb]{0.501961,0,0}{\textbf{224.9343}}  & \textcolor[rgb]{0.501961,0,0}{\textbf{22.7478}}  & \textbf{46.0225}  & \textbf{53.7072}\tabularnewline
\hline
\textcolor[rgb]{0.501961,0,0}{\textbf{458.2226}}  & \textcolor[rgb]{0.501961,0,0}{\textbf{305.6334}}  & \textcolor[rgb]{0.501961,0,0}{\textbf{214.9755}}  & \textcolor[rgb]{0.501961,0,0}{\textbf{21.6373}}  & \textbf{45.6650}  & \textbf{26.3564}\tabularnewline
\hline
\textcolor[rgb]{0.501961,0,0}{\textbf{435.2661}}  & \textcolor[rgb]{0.501961,0,0}{\textbf{13.1683}}  & \textcolor[rgb]{0.501961,0,0}{\textbf{45.7319}}  & \textcolor[rgb]{0.501961,0,0}{\textbf{21.3150}}  & \textbf{25.4884}  & \textbf{25.6006}\tabularnewline
\hline
\textcolor[rgb]{0.501961,0,0}{\textbf{392.6272}}  & \textcolor[rgb]{0.501961,0,0}{\textbf{11.1482}}  & \textcolor[rgb]{0.501961,0,0}{\textbf{17.8374}}  & \textcolor[rgb]{0.501961,0,0}{\textbf{19.2596}}  & \textbf{18.0820}  & \textbf{19.3930}\tabularnewline
\hline
3.6928  & \textcolor[rgb]{0.501961,0,0}{\textbf{9.0674}}  & \textcolor[rgb]{0.501961,0,0}{\textbf{11.7423}}  & \textcolor[rgb]{0.501961,0,0}{\textbf{10.3580}}  & \textbf{15.5876}  & \textbf{14.6107}\tabularnewline
\hline
3.5809  & \textcolor[rgb]{0.501961,0,0}{\textbf{7.9537}}  & \textcolor[rgb]{0.501961,0,0}{\textbf{8.0837}}  & \textcolor[rgb]{0.501961,0,0}{\textbf{9.9668}}  & \textbf{12.9088}  & \textbf{12.0980}\tabularnewline
\hline
3.5449  & 1.7375  &  1.7988 & \textcolor[rgb]{0.501961,0,0}{\textbf{9.5028}}  & \textbf{10.5568}  & \textbf{8.5791}\tabularnewline
\hline
3.4579  & 1.7326  & 1.7895  & \textcolor[rgb]{0.501961,0,0}{\textbf{8.5483}}  & \textbf{4.7840}  & \textbf{7.0596}\tabularnewline
\hline
3.4312  & 1.7015  &  1.7388 & \textcolor[rgb]{0.501961,0,0}{\textbf{5.5931}} & \textbf{4.3896}  & \textbf{4.5411}\tabularnewline
\hline
3.3829  & 1.6957  & 1.7242  & \textcolor[rgb]{0.501961,0,0}{\textbf{3.5968}}  & \textbf{4.3843}  & \textbf{3.6744}\tabularnewline
\hline
3.3701  & 1.6751  & 1.6724 & 1.8215 & 1.7944  & 1.7468\tabularnewline
\hline
\end{tabular}\caption{Absolute values of the largest eigenvalues of the empirical covariance
matrix at various lags with parameters: $n=450$, $T=500$, $k=2$,
$q=2$, $\beta=1.0$, $\sigma_{\mathit{f}}^{2}$ = 4, $\sigma^2=1$
and $\sigma_{\ep}^{2}=0.25$. Note that $c=0.9$, $b_c=(1+\sqrt c)^2=3.7974$ and $d_c=1.8573$.
When $\sigma^{2}=1$ is unknown, one has $\hat\sigma^2=0.9894$, $\hat b_c=3.8851$ and $\hat d_c=1.8616$.}
\end{table}
\appendix
\section{Some proofs}
\setcounter{equation}{0}
\def\theequation{\thesection.\arabic{equation}}
\setcounter{subsection}{0}
\subsection{Proof of Lemma \ref{Cp}}
Define $\bbM_k=\bbM-\bgma_k(\bgma_{k+\tau}+\bgma_{k-\tau})^*-(\bgma_{k+\tau}+\bgma_{k-\tau})\bgma_k^*$, and
\bqn
\bbM_{k,k+\tau,\cdots,k+l\tau}=\bbM_{k,k+\tau,\cdots,k+(l-1)\tau}-\bgma_{k+(l+1)\tau}\bgma_{k+l\tau}^*-\bgma_{k+l\tau}\bgma_{k+(l+1)\tau}^*, l\geq 1.
\eqn
Suppose that $i\ge j$, then we have
\bqn
&&\bgma_i^*\bbB_2\Big(\ell\bbI-\bbB_2^*\bbM\bbB_2\Big)^{-1}\bbB_2^*\bgma_j\\
&=&\bgma_i^*\bbB_2\Big(\ell\bbI-\bbB_2^*\bbM_j\bbB_2-\bbB_2^*(\bgma_{j+\tau}+\bgma_{j-\tau})\bgma_j^*\bbB_2
-\bbB_2^*\bgma_j(\bgma_{j+\tau}+\bgma_{j-\tau}^*)\bbB_2\Big)^{-1}\bbB_2^*\bgma_j\\
&=&\frac{\bgma_i^*\bbB_2\Big(\ell\bbI-\bbB_2^*\bbM_j\bbB_2-\bbB_2^*(\bgma_{j+\tau}+\bgma_{j-\tau})\bgma_j^*\bbB_2\Big)^{-1}\bbB_2^*\bgma_j}
{1-(\bgma_{j+\tau}+\bgma_{j-\tau}^*)\bbB_2\Big(\ell\bbI-\bbB_2^*\bbM_j\bbB_2-\bbB_2^*(\bgma_{j+\tau}+\bgma_{j-\tau})\bgma_j^*\bbB_2\Big)^{-1}\bbB_2^*\bgma_j}\\
&=&\frac{\bgma_i^*\bbB_2\Big((\ell\bbI-\bbB_2^*\bbM_j\bbB_2)^{-1}+\frac{(\ell\bbI-\bbB_2^*\bbM_j\bbB_2)^{-1}\bbB_2^*(\bgma_{j+\tau}+\bgma_{j-\tau})\bgma_j^*\bbB_2(\ell\bbI-\bbB_2^*\bbM_j\bbB_2)^{-1}}{1-\bgma_j^*\bbB_2(\ell\bbI-\bbB_2^*\bbM_j\bbB_2)^{-1}\bbB_2^*(\bgma_{j+\tau}+\bgma_{j-\tau})}\Big)\bbB_2^*\bgma_j}
{1-(\bgma_{j+\tau}+\bgma_{j-\tau}^*)\bbB_2\Big((\ell\bbI-\bbB_2^*\bbM_j\bbB_2)^{-1}+\frac{(\ell\bbI-\bbB_2^*\bbM_j\bbB_2)^{-1}\bbB_2^*(\bgma_{j+\tau}+\bgma_{j-\tau})\bgma_j^*\bbB_2(\ell\bbI-\bbB_2^*\bbM_j\bbB_2)^{-1}}{1-\bgma_j^*\bbB_2(\ell\bbI-\bbB_2^*\bbM_j\bbB_2)^{-1}\bbB_2^*(\bgma_{j+\tau}+\bgma_{j-\tau})}\Big)\bbB_2^*\bgma_j}
\eqn
\bqn
&=\frac{\bgma_i^*\bbB_2\Big((\ell\bbI-\bbB_2^*\bbM_j\bbB_2)^{-1}+(\ell\bbI-\bbB_2^*\bbM_j\bbB_2)^{-1}\bbB_2^*(\bgma_{j+\tau}+\bgma_{j-\tau})\bgma_j^*\bbB_2(\ell\bbI-\bbB_2^*\bbM_j\bbB_2)^{-1}\Big)\bbB_2^*\bgma_j}
{1-(\bgma_{j+\tau}+\bgma_{j-\tau}^*)\bbB_2\Big((\ell\bbI-\bbB_2^*\bbM_j\bbB_2)^{-1}+(\ell\bbI-\bbB_2^*\bbM_j\bbB_2)^{-1}\bbB_2^*(\bgma_{j+\tau}+\bgma_{j-\tau})\bgma_j^*\bbB_2(\ell\bbI-\bbB_2^*\bbM_j\bbB_2)^{-1}\Big)\bbB_2^*\bgma_j}+o_{a.s.}(1)&
\eqn
\bqn
=\left\{\begin{array}{l}\frac{-\frac{cm}2}
{1-\frac{c^2m^2}{2x_1}}+o_{a.s.}(1), \quad i=j\\
\frac{-\frac{cm}2}{1-\frac{c^2m^2}{2x_1}}\bgma_i^*\bbB_2(\ell\bbI-\bbB_2^*\bbM_j\bbB_2)^{-1}\bbB_2^*(\bgma_{j+\tau}+\bgma_{j-\tau})
+o_{a.s.}(1), \quad \mbox{otherwise}.
\end{array}\right.
\eqn
Next, we have
\bqn
&&\bgma_i^*\bbB_2(\ell\bbI-\bbB_2^*\bbM_j\bbB_2)^{-1}\bbB_2^*\bgma_{j+\tau}\\
&=&\bgma_i^*\bbB_2(\ell\bbI-\bbB_2^*\bbM_{j,j+\tau}\bbB_2-\bbB_2^*\bgma_{j+2\tau}\bgma_{j+\tau}^*\bbB_2-\bbB_2^*\bgma_{j+\tau}\bgma_{j+2\tau}^*\bbB_2)^{-1}\bbB_2^*\bgma_{j+\tau}\\
&=&\frac{\bgma_i^*\bbB_2(\ell\bbI-\bbB_2^*\bbM_{j,j+\tau}\bbB_2-\bbB_2^*\bgma_{j+2\tau}\bgma_{j+\tau}^*\bbB_2)^{-1}\bbB_2^*\bgma_{j+\tau}}
{1-\bgma_{j+2\tau}^*\bbB_2(\ell\bbI-\bbB_2^*\bbM_{j,j+\tau}\bbB_2-\bbB_2^*\bgma_{j+2\tau}\bgma_{j+\tau}^*\bbB_2)^{-1}\bbB_2^*\bgma_{j+\tau}}\\
&=&\frac{\bgma_i^*\bbB_2\Big((\ell\bbI-\bbB_2^*\bbM_{j,j+\tau}\bbB_2)^{-1}+\frac{(\ell\bbI-\bbB_2^*\bbM_{j,j+\tau}\bbB_2)^{-1}\bbB_2^*\bgma_{j+2\tau}\bgma_{j+\tau}^*\bbB_2(\ell\bbI-\bbB_2^*\bbM_{j,j+\tau}\bbB_2)^{-1}}{1-\bgma_{j+\tau}^*\bbB_2(\ell\bbI-\bbB_2^*\bbM_{j,j+\tau}\bbB_2)^{-1}\bbB_2^*\bgma_{j+2\tau}}\Big)\bbB_2^*\bgma_{j+\tau}}
{1-\bgma_{j+2\tau}^*\bbB_2\Big((\ell\bbI-\bbB_2^*\bbM_{j,j+\tau}\bbB_2)^{-1}+\frac{(\ell\bbI-\bbB_2^*\bbM_{j,j+\tau}\bbB_2)^{-1}\bbB_2^*\bgma_{j+2\tau}\bgma_{j+\tau}^*\bbB_2(\ell\bbI-\bbB_2^*\bbM_{j,j+\tau}\bbB_2)^{-1}}{1-\bgma_{j+\tau}^*\bbB_2(\ell\bbI-\bbB_2^*\bbM_{j,j+\tau}\bbB_2)^{-1}\bbB_2^*\bgma_{j+2\tau}}\Big)\bbB_2^*\bgma_{j+\tau}}
\eqn
\bqn
&=\frac{\bgma_i^*\bbB_2\Big((\ell\bbI-\bbB_2^*\bbM_{j,j+\tau}\bbB_2)^{-1}+(\ell\bbI-\bbB_2^*\bbM_{j,j+\tau}\bbB_2)^{-1}\bbB_2^*\bgma_{j+2\tau}\bgma_{j+\tau}^*\bbB_2(\ell\bbI-\bbB_2^*\bbM_{j,j+\tau}\bbB_2)^{-1}\Big)\bbB_2^*\bgma_{j+\tau}}
{1-\bgma_{j+2\tau}^*\bbB_2\Big((\ell\bbI-\bbB_2^*\bbM_{j,j+\tau}\bbB_2)^{-1}+(\ell\bbI-\bbB_2^*\bbM_{j,j+\tau}\bbB_2)^{-1}\bbB_2^*\bgma_{j+2\tau}\bgma_{j+\tau}^*\bbB_2(\ell\bbI-\bbB_2^*\bbM_{j,j+\tau}\bbB_2)^{-1}\Big)\bbB_2^*\bgma_{j+\tau}}+o_{a.s.}(1)&
\eqn
\bqn
=\left\{\begin{array}{l}\frac{-\frac{cm}2}
{1-\frac{c^2m^2}{4x_1}}+o_{a.s.}(1), \quad i=j+\tau\\
\frac{-\frac{cm}2}{1-\frac{c^2m^2}{4x_1}}\bgma_i^*\bbB_2(\ell\bbI-\bbB_2^*\bbM_{j,j+\tau}\bbB_2)^{-1}\bbB_2^*\bgma_{j+2\tau}+o_{a.s.}(1), \quad \mbox{otherwise}.
\end{array}\right.
\eqn
\bqn
&&\bgma_i^*\bbB_2(\ell\bbI-\bbB_2^*\bbM_j\bbB_2)^{-1}\bbB_2^*\bgma_{j-\tau}\\
&=&\bgma_i^*\bbB_2(\ell\bbI-\bbB_2^*\bbM_{j,j-\tau}\bbB_2-\bbB_2^*\bgma_{j-2\tau}\bgma_{j-\tau}^*\bbB_2-\bbB_2^*\bgma_{j-\tau}\bgma_{j-2\tau}^*\bbB_2)^{-1}\bbB_2^*\bgma_{j-\tau}\\
&=&\frac{\bgma_i^*\bbB_2(\ell\bbI-\bbB_2^*\bbM_{j,j-\tau}\bbB_2-\bbB_2^*\bgma_{j-2\tau}\bgma_{j-\tau}^*\bbB_2)^{-1}\bbB_2^*\bgma_{j-\tau}}
{1-\bgma_{j-2\tau}^*\bbB_2(\ell\bbI-\bbB_2^*\bbM_{j,j-\tau}\bbB_2-\bbB_2^*\bgma_{j-2\tau}\bgma_{j-\tau}^*\bbB_2)^{-1}\bbB_2^*\bgma_{j-\tau}}\\
&=&\frac{\bgma_i^*\bbB_2\Big((\ell\bbI-\bbB_2^*\bbM_{j,j-\tau}\bbB_2)^{-1}+\frac{(\ell\bbI-\bbB_2^*\bbM_{j,j-\tau}\bbB_2)^{-1}\bbB_2^*\bgma_{j-2\tau}\bgma_{j-\tau}^*\bbB_2(\ell\bbI-\bbB_2^*\bbM_{j,j-\tau}\bbB_2)^{-1}}{1-\bgma_{j-\tau}^*\bbB_2(\ell\bbI-\bbB_2^*\bbM_{j,j-\tau}\bbB_2)^{-1}\bbB_2^*\bgma_{j-2\tau}}\Big)\bbB_2^*\bgma_{j-\tau}}
{1-\bgma_{j-2\tau}^*\bbB_2\Big((\ell\bbI-\bbB_2^*\bbM_{j,j-\tau}\bbB_2)^{-1}+\frac{(\ell\bbI-\bbB_2^*\bbM_{j,j-\tau}\bbB_2)^{-1}\bbB_2^*\bgma_{j-2\tau}\bgma_{j-\tau}^*\bbB_2(\ell\bbI-\bbB_2^*\bbM_{j,j-\tau}\bbB_2)^{-1}}{1-\bgma_{j-\tau}^*\bbB_2(\ell\bbI-\bbB_2^*\bbM_{j,j-\tau}\bbB_2)^{-1}\bbB_2^*\bgma_{j-2\tau}}\Big)\bbB_2^*\bgma_{j-\tau}}
\eqn
\bqn
&=\frac{\bgma_i^*\bbB_2\Big((\ell\bbI-\bbB_2^*\bbM_{j,j-\tau}\bbB_2)^{-1}+(\ell\bbI-\bbB_2^*\bbM_{j,j-\tau}\bbB_2)^{-1}\bbB_2^*\bgma_{j-2\tau}\bgma_{j-\tau}^*\bbB_2(\ell\bbI-\bbB_2^*\bbM_{j,j-\tau}\bbB_2)^{-1}\Big)\bbB_2^*\bgma_{j-\tau}}
{1-\bgma_{j-2\tau}^*\bbB_2\Big((\ell\bbI-\bbB_2^*\bbM_{j,j-\tau}\bbB_2)^{-1}+(\ell\bbI-\bbB_2^*\bbM_{j,j-\tau}\bbB_2)^{-1}\bbB_2^*\bgma_{j-2\tau}\bgma_{j-\tau}^*\bbB_2(\ell\bbI-\bbB_2^*\bbM_{j,j-\tau}\bbB_2)^{-1}\Big)\bbB_2^*\bgma_{j-\tau}}+o_{a.s.}(1)&
\eqn
\bqn
=\frac{-\frac{cm}2}{1-\frac{c^2m^2}{4x_1}}\bgma_i^*\bbB_2(\ell\bbI-\bbB_2^*\bbM_{j,j-\tau}\bbB_2)^{-1}\bbB_2^*\bgma_{j-2\tau}+o_{a.s.}(1).
\eqn
Note that $\left|\frac{-\frac{cm}2}{1-\frac{c^2m^2}{4x_1}}\right|=\left|\frac{a}{x_1}\right|<1$, by induction, we have $$\bgma_i^*\bbB_2(\ell\bbI-\bbB_2^*\bbM_j\bbB_2)^{-1}\bbB_2^*\bgma_{j-\tau}=o_{a.s.}(1).$$
Then result then follows by induction. By symmetry, it holds when $i<j$. The proof of the lemma is complete.
\subsection{Proof of Lemma \ref{gAx}}
Let $\bbA_k^{-1}\bbx=\bbb=(b_1,\cdots,b_n)'$. Noting $|\ep_{ij}|<C$ and $\rE|\ep_{ij}|^2=1$,
we have
\bqn
&&\rE(\bgma_k^*\bbA_k^{-1}\bbx)^{2r}
=\frac{1}{2^r T^r}\rE(\sum_{i=1}^n \ep_{ki}b_i)^{2r}\\
&\leq &\frac{1}{2^r T^r}\rE\sum_{l=1}^{r}\sum_{1\leq j_1< \cdots< j_l\leq n}\sum_{i_1+\cdots+i_l=2r} \frac{(2r)!}{i_1 !\cdots i_l !} \ep_{kj_1}^{i_1}b_{j_1}^{i_1}\cdots\ep_{kj_l}^{i_l}b_{j_1}^{i_l}\\
&=&\frac{1}{2^r T^r}\rE\sum_{l=1}^{r}\sum_{1\leq j_1< \cdots< j_l\leq n}\sum_{\begin{subarray}{l}
i_1+\cdots+i_l=2r,\\i_1\geq2,\cdots,i_l\geq 2\end{subarray}} \frac{(2r)!}{i_1 !\cdots i_l !} \ep_{kj_1}^{i_1}b_{j_1}^{i_1}\cdots\ep_{kj_l}^{i_l}b_{j_1}^{i_l}\\
&\le&\frac{K}{2^r T^r} \rE\sum_{l=1}^{r}\sum_{1\leq j_1< \cdots<j_l\le n}\sum_{\begin{subarray}{l}
i_1+\cdots+i_l=2r,\\i_1\geq2,\cdots,i_l\geq 2\end{subarray}} \frac{(2r)!}{i_1 !\cdots i_l !} |b_{j_1}|^{i_1}\cdots |b_{j_l}|^{i_l}.
\eqn
By $\sum_{j=1}^n |b_j|^2=\|\bbA_k^{-1}\bbx\|^2$ and Cauchy-Schwartz inequality, we have
\bqn
&&\sum_{1\leq j_1< \cdots<j_l\leq n}\sum_{\begin{subarray}{l}
i_1+\cdots+i_l=2r,\\i_1\geq2,\cdots,i_l\geq 2\end{subarray}} \frac{(2r)!}{i_1 !\cdots i_l !} |b_{j_1}|^{i_1}\cdots |b_{j_l}|^{i_l}\\
&\leq& \sum_{\begin{subarray}{l}
i_1+\cdots+i_l=2r,\\i_1\geq2,\cdots,i_l\geq 2\end{subarray}} \frac{(2r)!}{i_1 !\cdots i_l ! }(\sum_{j=1}^n |b_j|^2)^r\\
&\leq&l^{2r}||\bbA_k^{-1}||_o^{2r} ||\bbx||^{2r}\\
&\leq& \frac{l^{2r}}{\eta^{2r}}.
\eqn
Here $\eta:=\ell-d_c>0$.
Therefore, we have \bqn
\rE|\bgma_k^*\bbA_k^{-1}\bbx|^{2r}\le KT^{-r}
\eqn
for some $K>0$. The proof of the lemma is complete.
\subsection{Proof of Lemma \ref{xAy}}
\proof. Let $\bbx,\bby\in \mathbb{C}_1^n$ be given. Define $\bbA_{k,k+\tau}=\bbA_k-\bgma_{k+\tau}\bgma_{k+2\tau}^*-\bgma_{k+2\tau}\bgma_{k+\tau}^*$,\\
$\tilde\bbA_k=\bbA_k+\bgma_k(\bgma_{k+\tau}+\bgma_{k-\tau})^*$ and $\tilde\bbA_{k,k+\tau}=\bbA_k-\bgma_{k+2\tau}\bgma_{k+\tau}^*$. First we have
\bqn
&&\bbx^*\bbA^{-1}\bby-\rE\bbx^*\bbA^{-1}\bby\\
&=&\bbx^*\sum_{k=1}^T(\rE_k-\rE_{k-1})\big(\bbA^{-1}-\bbA_k^{-1}(\ell)\big)\bby\\
&=&\sum_{k=1}^T(\rE_k-\rE_{k-1})\big(-\frac{\bbx^*\tilde\bbA_k^{-1}(\bgma_{k+\tau}+\bgma_{k-\tau})\bgma_k^*\tilde\bbA_k^{-1}\bby}{1+\bgma_k^*\tilde\bbA_k^{-1}(\bgma_{k+\tau}+\bgma_{k-\tau})}
-\frac{\bbx^*\bbA_k^{-1}\bgma_k(\bgma_{k+\tau}+\bgma_{k-\tau})^*\bbA_k^{-1}\bby}{1+(\bgma_{k+\tau}+\bgma_{k-\tau})^*\bbA_k^{-1}\bgma_k}\big)\\
&\equiv&\sum_{k=1}^T(\rE_k-\rE_{k-1})\big(-\alpha_{k1}-\alpha_{k2}).
\eqn
Using Lemma \ref{Burkholder}, we have, for $i=1,2$
\bqa\label{bur}
&&\rE|\sum_{k=1}^T(\rE_k-\rE_{k-1})\alpha_{ki}|^{2l}\non
&\le& K_l \left[\rE\left(\sum_{k=1}^{T}\rE_{k-1}|(\rE_k-\rE_{k-1})\alpha_{ki}|^2\right)^l+\sum_{k=1}^{T}\rE|(\rE_k-\rE_{k-1})\alpha_{ki}|^{2l}\right]\non
&\le&K'_l \left[\rE\left(\sum_{k=1}^{T}\rE_{k-1}|\rE_k\alpha_{ki}|^2+\sum_{k=1}^{T}\rE_{k-1}|\alpha_{ki}|^2\right)^l+\sum_{k=1}^{T}\rE|\rE_k\alpha_{ki}|^{2l}+\sum_{k=1}^{T}\rE|\rE_{k-1}\alpha_{ki}|^{2l}\right]\non
&\le&2^{l}K'_l
\left[\rE\left(\sum_{k=1}^{T}\rE_{k-1}|\alpha_{ki}|^2\right)^l+\sum_{k=1}^{T}\rE|\alpha_{ki}|^{2l}\right].
\eqa
Note that
\bqn
\bbA_k^{-1}=(\tilde\bbA_{k,k+\tau}+\bgma_{k+2\tau}\bgma_{k+\tau}^*)^{-1}=\tilde\bbA_{k,k+\tau}^{-1}-\frac{\tilde\bbA_{k,k+\tau}^{-1}\bgma_{k+2\tau}\bgma_{k+\tau}^*\tilde\bbA_{k,k+\tau}^{-1}}{1+\bgma_{k+\tau}^*\tilde\bbA_{k,k+\tau}^{-1}\bgma_{k+2\tau}}.
\eqn
Hence, we have
\bqn
\bgma_{k+\tau}^*\bbA_k^{-1}=\bgma_{k+\tau}^*\tilde\bbA_{k,k+\tau}^{-1}-\frac{\bgma_{k+\tau}^*\tilde\bbA_{k,k+\tau}^{-1}\bgma_{k+2\tau}\bgma_{k+\tau}^*\tilde\bbA_{k,k+\tau}^{-1}}{1+\bgma_{k+\tau}^*\tilde\bbA_{k,k+\tau}^{-1}\bgma_{k+2\tau}}
=\frac{\bgma_{k+\tau}^*\tilde\bbA_{k,k+\tau}^{-1}}{1+\bgma_{k+\tau}^*\tilde\bbA_{k,k+\tau}^{-1}\bgma_{k+2\tau}}.
\eqn
Next, we have
\bqn
\bgma_{k+\tau}^*\tilde\bbA_{k,k+\tau}^{-1}&=&\bgma_{k+\tau}^*\bbA_{k,k+\tau}^{-1}-\frac{\bgma_{k+\tau}^*\bbA_{k,k+\tau}^{-1}\bgma_{k+\tau}\bgma_{k+2\tau}^*\bbA_{k,k+\tau}^{-1}}{1+\bgma_{k+2\tau}^*\bbA_{k,k+\tau}^{-1}\bgma_{k+\tau}}\\
&=&\bgma_{k+\tau}^*\bbA_{k,k+\tau}^{-1}-\frac{cm}{2}\bgma_{k+2\tau}^*\bbA_{k,k+\tau}^{-1}+R_{k1},
\eqn
where
\bqn
R_{k1}&=&\frac{cm}{2}\bgma_{k+2\tau}^*\bbA_{k,k+\tau}^{-1}-\frac{\bgma_{k+\tau}^*\bbA_{k,k+\tau}^{-1}\bgma_{k+\tau}\bgma_{k+2\tau}^*\bbA_{k,k+\tau}^{-1}}{1+\bgma_{k+2\tau}^*\bbA_{k,k+\tau}^{-1}\bgma_{k+\tau}}\\
&=&\left(\frac{\frac{cm}{2}-\bgma_{k+\tau}^*\bbA_{k,k+\tau}^{-1}\bgma_{k+\tau}+\frac{cm}{2}\bgma_{k+2\tau}^*\bbA_{k,k+\tau}^{-1}\bgma_{k+\tau}}{1+\bgma_{k+2\tau}^*\bbA_{k,k+\tau}^{-1}\bgma_{k+\tau}}\right)\bgma_{k+2\tau}^*\bbA_{k,k+\tau}^{-1}.
\eqn
Substitute back, we obtain
\bqa\label{eq2}
&&\bgma_{k+\tau}^*\bbA_k^{-1}\bby\non
&=&\frac{\bgma_{k+\tau}^*\bbA_{k,k+\tau}^{-1}\bby-\frac{\bgma_{k+\tau}^*\bbA_{k,k+\tau}^{-1}\bgma_{k+\tau}}{1+\bgma_{k+2\tau}^*\bbA_{k,k+\tau}^{-1}\bgma_{k+\tau}}\bgma_{k+2\tau}^*\bbA_{k,k+\tau}^{-1}\bby}{1+\bgma_{k+\tau}^*\bbA_{k,k+\tau}^{-1}\bgma_{k+2\tau}-\frac{cm}{2}\bgma_{k+2\tau}^*\bbA_{k,k+\tau}^{-1}\bgma_{k+2\tau}+R_{k1}\bgma_{k+2\tau}},
\eqa
with
\bqn
R_{k1}\bgma_{k+2\tau}=\left(\frac{\frac{cm}{2}-\bgma_{k+\tau}^*\bbA_{k,k+\tau}^{-1}\bgma_{k+\tau}+\frac{cm}{2}\bgma_{k+2\tau}^*\bbA_{k,k+\tau}^{-1}\bgma_{k+\tau}}{1+\bgma_{k+2\tau}^*\bbA_{k,k+\tau}^{-1}\bgma_{k+\tau}}\right)\bgma_{k+2\tau}^*\bbA_{k,k+\tau}^{-1}\bgma_{k+2\tau}.
\eqn
Note that
\bqn
&&\lim_{n\to\infty}\left|\frac{\bgma_{k+\tau}^*\bbA_{k,k+\tau}^{-1}\bgma_{k+\tau}}{\big(1+\bgma_{k+2\tau}^*\bbA_{k,k+\tau}^{-1}\bgma_{k+\tau}\big)\big(1+\bgma_{k+\tau}^*\bbA_{k,k+\tau}^{-1}\bgma_{k+2\tau}-\frac{cm}{2}\bgma_{k+2\tau}^*\bbA_{k,k+\tau}^{-1}\bgma_{k+2\tau}+R_{k1}\bgma_{k+2\tau}\big)}\right|\\
&=&\left|\frac{\frac{cm}{2}}{1-\frac{cm}{2}\frac{cm}{2x_1}}\right|=\left|\frac{\frac{cm}2}{x_1}\right|<1
\eqn
and
\bqn
\lim_{n\to\infty}1+\bgma_{k+\tau}^*\bbA_{k,k+\tau}^{-1}\bgma_{k+2\tau}-\frac{cm}{2}\bgma_{k+2\tau}^*\bbA_{k,k+\tau}^{-1}\bgma_{k+2\tau}+R_{k1}\bgma_{k+2\tau}
=1-\frac{cm}{2}\frac{cm}{2x_1},
\eqn
which is bounded. Using induction, we have $|\bgma_{k+\tau}^*\bbA_k^{-1}\bby|\le K|\bgma_{k+\tau}^*\bbA_{k,k+\tau}^{-1}\bby|$.\\
Similarly, we have $|\bgma_{k-\tau}^*\bbA_k^{-1}\bby|\le K|\bgma_{k-\tau}^*\bbA_{k,k-\tau}^{-1}\bby|$,
$|\bbx^*\bbA_k^{-1}\bgma_{k+\tau}|\le K|\bbx^*\bbA_{k,k+\tau}^{-1}\bgma_{k+\tau}|$ and\\
$|\bbx^*\bbA_k^{-1}\bgma_{k-\tau}|\le K|\bbx^*\bbA_{k,k-\tau}^{-1}\bgma_{k-\tau}|$.\\
Therefore, by noting $(\bgma_{k+\tau}+\bgma_{k-\tau})^*\bbA_k^{-1}\bgma_k=o_{a.s.}(1)$,
$|\ep_{it}|<C$, $\rE|\ep_{it}|^2=1$ and $\bbx^*\bbA_k^{-1}\bar\bbA_k^{-1}\bbx$ being bounded, we have
\bqn
&&\rE\left(\sum_{k=1}^{T}\rE_{k-1}|\alpha_{k2}|^2\right)^l\\
&\le&K\rE\left(\sum_{k=1}^{T}\rE_{k-1}|\bbx^*\bbA_k^{-1}\bgma_k\bgma_{k+\tau}\bbA_k^{-1}\bby|^2\right)^l\\
&=&K\rE\left(\sum_{k=1}^{T}\frac1{2T}\rE_{k-1}\bbx^*\bbA_k^{-1}\bar\bbA_k^{-1}\bbx|\bgma_{k+\tau}\bbA_k^{-1}\bby|^2\right)^l\\
&\le&K\max_k\rE|\bgma_{k+\tau}\bbA_{k,k+\tau}^{-1}\bby|^{2l}\\
&\le&\frac K {T^l}
\eqn
and
\bqn
&&\sum_{k=1}^{T}\rE|\alpha_{k2}|^{2l}\\
&=&\sum_{k=1}^{T}\rE\left(|\bbx^*\bbA_k^{-1}\bgma_k\bgma_{k+\tau}\bbA_k^{-1}\bby|^2\right)^l\\
&\le&\frac K {T^{l-1}}\max_k\rE|\bgma_{k+\tau}\bbA_{k,k+\tau}^{-1}\bby|^{2l}\\
&\le&\frac K {T^{2l-1}}
\eqn
For $i=1$, by $\tilde\bbA_k^{-1}=\bbA_k^{-1}-\frac{\bbA_k^{-1}\bgma_k(\bgma_{k+\tau}+\bgma_{k-\tau})^*\bbA_k^{-1}}{1+(\bgma_{k+\tau}+\bgma_{k-\tau})^*\bbA_k^{-1}\bgma_k}$,
we have
\bqn
&&\bbx^*\tilde\bbA_k^{-1}(\bgma_{k+\tau}+\bgma_{k-\tau})\bgma_k^*\tilde\bbA_k^{-1}\bby\\
&=&\bbx^*\Big(\bbA_k^{-1}-\frac{\bbA_k^{-1}\bgma_k(\bgma_{k+\tau}+\bgma_{k-\tau})^*\bbA_k^{-1}}{1+(\bgma_{k+\tau}+\bgma_{k-\tau})^*\bbA_k^{-1}\bgma_k}\Big)
(\bgma_{k+\tau}+\bgma_{k-\tau})\bgma_k^*\Big(\bbA_k^{-1}-\frac{\bbA_k^{-1}\bgma_k(\bgma_{k+\tau}+\bgma_{k-\tau})^*\bbA_k^{-1}}{1+(\bgma_{k+\tau}+\bgma_{k-\tau})^*\bbA_k^{-1}\bgma_k}\Big)\bby\\
&=&\bbx^*\bbA_k^{-1}(\bgma_{k+\tau}+\bgma_{k-\tau})\bgma_k^*\bbA_k^{-1}\bby\\
&&-\frac{(\bgma_{k+\tau}+\bgma_{k-\tau})^*\bbA_k^{-1}
(\bgma_{k+\tau}+\bgma_{k-\tau})}{1+(\bgma_{k+\tau}+\bgma_{k-\tau})^*\bbA_k^{-1}\bgma_k}\bbx^*\bbA_k^{-1}\bgma_k\bgma_k^*\bbA_k^{-1}\bby\\
&&-\frac{\bgma_k^*
\bbA_k^{-1}\bgma_k}{1+(\bgma_{k+\tau}+\bgma_{k-\tau})^*\bbA_k^{-1}\bgma_k}\bbx^*\bbA_k^{-1}(\bgma_{k+\tau}+\bgma_{k-\tau})(\bgma_{k+\tau}+\bgma_{k-\tau})^*\bbA_k^{-1}\bby\\
&&+\frac{\bgma_k^*\bbA_k^{-1}\bgma_k(\bgma_{k+\tau}+\bgma_{k-\tau})^*\bbA_k^{-1}
(\bgma_{k+\tau}+\bgma_{k-\tau})}{\big(1+(\bgma_{k+\tau}+\bgma_{k-\tau})^*\bbA_k^{-1}\bgma_k\big)^2}\bbx^*\bbA_k^{-1}\bgma_k(\bgma_{k+\tau}+\bgma_{k-\tau})^*\bbA_k^{-1}\bby\\
&=:&\alpha_{k11}-\alpha_{k12}-\alpha_{k13}+\alpha_{k14}
\eqn
It is easy to see that work on $\alpha_{k11}$ and $\alpha_{k14}$ is the same as that on $\alpha_{k_2}$.\\
For $\alpha_{k12}$, by Cauchy-Schwartz's inequality, we have
\bqn
&&\rE\left(\sum_{k=1}^{T}\rE_{k-1}|\alpha_{k12}|^2\right)^l\\
&\le&K\rE\left(\sum_{k=1}^{T}\rE_{k-1}|\bbx^*\bbA_k^{-1}\bgma_k|^2\rE_{k-1}|\bgma_k^*\bbA_k^{-1}\bby|^2\right)^l\\
&=&K\rE\left(\sum_{k=1}^{T}\frac1{4T^2}\rE_{k-1}\bbx^*\bbA_k^{-1}\bar\bbA_k^{-1}\bbx\rE_{k-1}\bby^*\bar\bbA_k^{-1}\bbA_k^{-1}\bby\right)^l\\
&\le&\frac K {T^l}
\eqn
and by $|\ep_{it}|<C$, $\rE|\ep_{it}|^2=1$ and $\bbx^*\bbA_k^{-1}\bar\bbA_k^{-1}\bbx$ being bounded, we have
\bqn
&&\sum_{k=1}^{T}\rE|\alpha_{k12}|^{2l}\\
&\le&K\sum_{k=1}^{T}\rE\left(|\bbx^*\bbA_k^{-1}\bgma_k|^2|\bgma_k^*\bbA_k^{-1}\bby|^2\right)^l\\
&\le&\frac K {T^{l-1}}\max_k\rE|\bgma_k\bbA_k^{-1}\bby|^{2l}\\
&\le&\frac K {T^{2l-1}}.
\eqn
By the fact that $|\bgma_{k\pm\tau}^*\bbA_k^{-1}\bby|\le K|\bgma_{k\pm\tau}^*\bbA_{k,k\pm\tau}^{-1}\bby|$ and
$|\bbx^*\bbA_k^{-1}\bgma_{k\pm\tau}|\le K|\bbx^*\bbA_{k,k\pm\tau}^{-1}\bgma_{k\pm\tau}|$, the similar result for $\alpha_{k13}$ follows by the same reason.\\
Substituting all the above results into (\ref{bur}) and choosing $l$ large enough, we have $\bbx^*\bbA^{-1}\bby-\rE\bbx^*\bbA^{-1}\bby\to0$ a.s.

Next, we want to show the convergence of $\rE\bbx^*\bbA^{-1}\bby$.\\
By
$$
\bbA=\sum_{k=1}^T(\bgma_k\bgma_{k+\tau}^*+\bgma_{k+\tau}\bgma_k^*)-\ell\bbI_n
$$
we have
$$
\bbI_n=\sum_{k=1}^T(\bgma_k\bgma_{k+\tau}^*\bbA^{-1}+\bgma_{k+\tau}\bgma_k^*\bbA^{-1})-\ell\bbA^{-1}.
$$
Multiplying $\bbx^*$ from left and $\bby$ from right and taking expectation, we obtain
\bqa
\bbx^*\bby&=&\sum_{k=1}^T(\rE\bbx^*\bgma_k\bgma_{k+\tau}^*\bbA^{-1}\bby+\rE\bbx^*\bgma_{k+\tau}\bgma_k^*\bbA^{-1}\bby)-\ell\rE\bbx^*\bbA^{-1}\bby\non
&=&\sum_{k=1}^T\rE\bbx^*(\bgma_{k+\tau}+\bgma_{k-\tau})\bgma_k^*\bbA^{-1}\bby-\ell\rE\bbx^*\bbA^{-1}\bby.
\label{eq1}
\eqa
By $\bbA=\tilde\bbA_k+(\bgma_{k+\tau}+\bgma_{k-\tau})\bgma_k^*$, $\tilde\bbA_k=\bbA_k+\bgma_k(\bgma_{k+\tau}+\bgma_{k-\tau})^*$,
$(\bgma_{k+\tau}+\bgma_{k-\tau})^*\bbA_k^{-1}\bgma_k=o_{a.s.}(1)$, $\bgma_k^*\bbA_k^{-1}\bgma_k=\frac{cm}2+o_{a.s.}(1)$, we have
\bqa
&&\bgma_k^*\bbA^{-1}\bby\non
&=&\frac{\bgma_k^*\tilde\bbA_k^{-1}\bby}{1+\bgma_k^*\tilde\bbA_k^{-1}(\bgma_{k+\tau}+\bgma_{k-\tau})}\non
&=&\frac{\bgma_k^*\Big(\bbA_k^{-1}-\frac{\bbA_k^{-1}\bgma_k(\bgma_{k+\tau}+\bgma_{k-\tau})^*\bbA_k^{-1}}{1+(\bgma_{k+\tau}+\bgma_{k-\tau})^*\bbA_k^{-1}\bgma_k}\Big)\bby}{1+\bgma_k^*\Big(\bbA_k^{-1}-\frac{\bbA_k^{-1}\bgma_k(\bgma_{k+\tau}+\bgma_{k-\tau})^*\bbA_k^{-1}}{1+(\bgma_{k+\tau}+\bgma_{k-\tau})^*\bbA_k^{-1}\bgma_k}\Big)(\bgma_{k+\tau}+\bgma_{k-\tau})}\non
&=&\frac{\bgma_k^*\bbA_k^{-1}\bby+\bgma_k^*\bbA_k^{-1}\bby(\bgma_{k+\tau}+\bgma_{k-\tau})^*\bbA_k^{-1}\bgma_k-\bgma_k^*\bbA_k^{-1}\bgma_k(\bgma_{k+\tau}+\bgma_{k-\tau})^*\bbA_k^{-1}\bby}
{[1+(\bgma_{k+\tau}+\bgma_{k-\tau})^*\bbA_k^{-1}\bgma_k][1+\bgma_k^*\bbA_k^{-1}(\bgma_{k+\tau}+\bgma_{k-\tau})]
-\bgma_k^*\bbA_k^{-1}\bgma_k(\bgma_{k+\tau}+\bgma_{k-\tau})^*\bbA_k^{-1}(\bgma_{k+\tau}+\bgma_{k-\tau})}\non
&=&\frac{-\frac{cm}2}{1-\frac{c^2m^2}{2x_1}}(\bgma_{k+\tau}+\bgma_{k-\tau})^*\bbA_k^{-1}\bby+o_{a.s.}(1).
\label{eq2}
\eqa
Next, we have
\bqn
&&\bgma_{k+\tau}^*\bbA_k^{-1}\bby\bbx^*\bgma_{k+\tau}\non
&=&\Big(\frac{(1+\bgma_{k+2\tau}^*\bbA_{k,k+\tau}^{-1}\bgma_{k+\tau})\bgma_{k+\tau}^*\bbA_{k,k+\tau}^{-1}\bby-\bgma_{k+\tau}^*\bbA_{k,k+\tau}^{-1}\bgma_{k+\tau}\bgma_{k+2\tau}^*\bbA_{k,k+\tau}^{-1}\bby}{(1+\bgma_{k+2\tau}^*\bbA_{k,k+\tau}^{-1}\bgma_{k+\tau})(1+\bgma_{k+\tau}^*\bbA_{k,k+\tau}^{-1}\bgma_{k+2\tau})-\bgma_{k+\tau}^*\bbA_{k,k+\tau}^{-1}\bgma_{k+\tau}\bgma_{k+2\tau}^*\bbA_{k,k+\tau}^{-1}\bgma_{k+2\tau}}\Big)\bbx^*\bgma_{k+\tau}\non
&=&\frac{\bgma_{k+\tau}^*\bbA_{k,k+\tau}^{-1}\bby\bbx^*\bgma_{k+\tau}+o_{a.s.}(1)}{1-\bgma_{k+\tau}^*\bbA_{k,k+\tau}^{-1}\bgma_{k+\tau}\bgma_{k+2\tau}^*\bbA_{k,k+\tau}^{-1}\bgma_{k+2\tau}}\non
&=&\frac{\bbx^*\bbA_{k,k+\tau}^{-1}\bby+o_{a.s.}(1)}{2T(1-\frac{c^2m^2}{4x_1})}.
\eqn
Similarly, we can show that $\bgma_{k-\tau}^*\bbA_k^{-1}\bby\bbx^*\bgma_{k-\tau}=\frac{\bbx^*\bbA_{k,k-\tau}^{-1}\bby+o_{a.s.}(1)}{2T(1-\frac{c^2m^2}{4x_1})}$,
$\bgma_{k+\tau}^*\bbA_k^{-1}\bby\bbx^*\bgma_{k-\tau}=o_{a.s.}(1)$ and $\bgma_{k-\tau}^*\bbA_k^{-1}\bby\bbx^*\bgma_{k+\tau}=o_{a.s.}(1)$.
Next, we will show that $\rE\bbx^*\bbA_{k,k\pm\tau}^{-1}\bby-\rE\bbx^*\bbA^{-1}\bby=o(1)$. By writing
$$\rE\bbx^*\bbA_{k,k\pm\tau}^{-1}\bby-\rE\bbx^*\bbA^{-1}\bby=\rE\bbx^*\bbA_{k,k\pm\tau}^{-1}\bby-\rE\bbx^*\bbA_{k}^{-1}\bby+\rE\bbx^*\bbA_{k}^{-1}\bby-\rE\bbx^*\bbA^{-1}\bby,$$ it is sufficient to show $\rE\bbx^*\bbA_k^{-1}\bby-\rE\bbx^*\bbA^{-1}\bby=o(1)$. Note that
\bqn
&&\rE\bbx^*\bbA_k^{-1}\bby-\rE\bbx^*\bbA^{-1}\bby\\
&=&\rE\big(\frac{\bbx^*\tilde\bbA_k^{-1}(\bgma_{k+\tau}+\bgma_{k-\tau})\bgma_k^*\tilde\bbA_k^{-1}\bby}{1+\bgma_k^*\tilde\bbA_k^{-1}(\bgma_{k+\tau}+\bgma_{k-\tau})}\big)
+\rE\big(\frac{\bbx^*\bbA_k^{-1}\bgma_k(\bgma_{k+\tau}+\bgma_{k-\tau})^*\bbA_k^{-1}\bby}{1+(\bgma_{k+\tau}+\bgma_{k-\tau})^*\bbA_k^{-1}\bgma_k}\big)\\
&=&\rE\alpha_{k1}+\rE\alpha_{k2}
\eqn
Previous calculation shows that $\rE|\alpha_{k1}|=o(1)$ and $\rE|\alpha_{k2}|=o(1)$.
Substituting these back to (\ref{eq1}) and (\ref{eq2}), we finish proving the lemma.
\bibliographystyle{plain}

\begin{thebibliography}{9}
 \bibitem{key1} Bai, Z.D., Liu H.X. and Wong, W.K. (2011)
 \newblock Asymptotic properties of eigenmatrices of a large sample covariance matrix.
 \newblock \emph{Ann. Appl. Probab.} 21, 1994--2015.
 \bibitem{key1} Bai, Z.D. and Silverstein, J.W.  (1998)
 \newblock No eigenvalues outside the support of the limiting spectral distribution of large-dimensional sample covariance matrices.
 \newblock \emph{Ann. Probab.} 26, 316--345.
\bibitem{key1} Bai, Z.D. and Silverstein, J.W. (2010)
 \newblock \emph{Spectral Analysis of Large Dimensional Random Matrices,} 2nd ed.
 \newblock Springer Verlag, New York.
 \bibitem{key1} Bai, Z.D. and Yao, J.F. (2008)
 \newblock Central limit theorems for eigenvalues in a spiked population model.
  \newblock {\it Ann. Inst. H. Poincar?????? Probab. Statist.} Volume 44, Number 3, 447--474.
\bibitem{key1} Bai, Z.D. and Wang, C. (2015)
 \newblock A note on the limiting spectral distribution of a symmetrized auto-cross covariance matrix.
 \newblock \emph{Statistics and Probability Letters}, 96, 333 -- 340.
  \bibitem{key1} Baik, Jinho and  Silverstein, J. W. (2006)
  \newblock Eigenvalues of large sample covariance matrices of spiked
 population models. \newblock {\it J. Multivariate Anal}. {97(6)}, 1382-1408.
  \bibitem{key1} Burkholder, D.L. (1973)
  \newblock Distribution function inequalities for martingales.
   \newblock \emph{Ann. Probab.} 1, 19--42.
 \bibitem{key1} Chung, K.L. (2001)
 \newblock \emph{A Course in Probability Theory,} 3rd ed.
 \newblock Academic Press, New York.
 \bibitem{key1} Johnstone, I. (2001)
 \newblock On the distribution of the largest eigenvalue in principal components analysis.
 \newblock \emph{Ann. Statist.} 29, 295 -- 327.
 \bibitem{key1} Jin, B. S., Wang, C., Bai, Z.D., Nair, K.K. and Harding, M.C. (2014)
 \newblock Limiting spectral distribution of a symmetrized auto-cross covariance matrix.
 \newblock \emph{Ann. Appl. Probab.} 24, 1199--1225.
 \bibitem{key1} Mar\v{c}enko, V.A. and Pastur, L.A. (1967)
 \newblock Distribution of eigenvalues for some sets of matrices.
 \newblock \emph{Mat. Sb.} 72, 507 -- 536.
 \bibitem{key1} Wang, C., Jin, B. S., Bai, Z.D., Nair, K.K. and Harding, M.C. (2015)
 \newblock Strong limit of the extreme eigenvalues of a symmetrized auto-cross covariance matrix.
 \newblock \emph{Ann. Appl. Probab.} 25, 3624 -- 3683.
\end{thebibliography}

\end{document}